\documentclass[12pt,leqno]{amsart}
\begin{document}
\bibliographystyle{plain}
\baselineskip16pt plus 1pt minus .5pt

\newtheorem{theorem}{Theorem}
\newtheorem{definition}[theorem]{Definition}
\newtheorem{convention}[theorem]{Convention}
\newtheorem{corollary}[theorem]{Corollary}
\newtheorem{observation}[theorem]{Observation}
\newtheorem{conjecture}[theorem]{Conjecture}
\newtheorem{odstavec}[theorem]{}
\newtheorem{example}[theorem]{Example}
\newtheorem{lemma}[theorem]{Lemma}
\newtheorem{sublemma}[theorem]{Sublemma}
\newtheorem{notation}[theorem]{Notation}
\newtheorem{proposition}[theorem]{Proposition}
\newtheorem{remark}[theorem]{Remark}
\newtheorem{remarks}[theorem]{Remarks}
\newtheorem{exercise}[theorem]{Exercise}
\newtheorem{claim}[theorem]{Claim}
\newtheorem{principle}[theorem]{Principle}

\def\boxtimes{\odot}
\def\invlim{{{\mathop{{\rm lim}}\limits_{\longleftarrow}}\hskip 1mm}}
\def\otm{{\hskip -.2mm \otimes \hskip -.2mm}}
\def\uXi{{\underline{\Xi}}}
\def\uGamma{{\underline{\Gamma}}}
\def\otexp#1#2{{#1}^{\otimes #2}}
\def\vlra{{\hbox{$-\hskip-1mm-\hskip-2mm\longrightarrow$}}}
\def\red{{\underline{\rm r}}}
\def\sred{{\underline{\mbox {\scriptsize r}}}}
\def\blue{{\underline{\rm b}}}
\def\sblue{{\mbox {\underline{\scriptsize b}}}}
\def\Vert{{\it Vert}}
\def\Leg{{\it Leg}}
\def\Flag{{\it Flag}}
\def\In{{\it In}}
\def\Out{{\it Out}}
\def\frprop{{\Gamma(\Xi)}}
\def\frpropDelta{{\Gamma_\Delta(\Xi)}}
\def\graph{{\sf G}}
\def\graphe{{\sf G}_e}
\def\Ass{\mbox{${\mathcal A}${\it ss}}}
\def\Comm{\mbox{${\mathcal C}${\hskip -.5mm \it om}}}
\def\zdvihatko{\rule{0pt}{.8em}}
\def\sB{{\rm s}{\sf B}}
\def\sfM{{\sf M}}\def\sfN{{\sf N}}
\def\DER{{\it Der}(\sfM,\End_V)}
\def\End{\hbox{${\mathcal E}\hskip -.1em {\it nd}$}}
\def\half{$\frac12$}
\def\uu{\underline}
\def\cric#1{\hbox{${}_{#1}\circ \hskip .2em$}}
\def\lra{\longrightarrow}
\def\cd{{\rm ({\bf c};{\bf d})}}
\def\GammaXi{{\Gamma(\Xi)}}
\def\genus{{{\rm gen}}} \def\pth{{\rm pth}} \def\grd{{\rm grd}}
\def\papert{\partial_{\it pert}}
\def\calD{{\mathcal D}}
\def\sfB{{\sf B}} \def\sfb{{\sf b}}  
\def\sfhb{{\mbox{$\frac 12${\sf b}}}} \def\hsfb{\sfhb}
\def\sshsfb{{\mbox{\scriptsize $\frac 12${\sf b}}}}
\def\calhB{{{\textstyle\frac12\mathcal B}}}
\def\calhBss{{{\frac12\mathcal B}}}
\def\BOX{\raisebox{1.2mm}%
         {\hskip .5mm $\fbox{\hphantom{\hglue .01mm}}$\hskip .5mm}}
\def\Box{\BOX}
\def\ssBox{\raisebox{1mm}%
         {\hskip .3mm $\fbox{\hphantom{\hglue -.2mm}}$\hskip .3mm}}
\def\sfhB{{\mbox{$\frac12$\sf B}}}
\def\sfp{{\sf p}} \def\sfs{{\sf s}}
\def\jcirc#1{{\hskip 1mm {}_{#1}\circ \hskip .2mm}}
\def\sfA{{\sf A}} \def\sfX{{\sf X}} \def\sfY{{\sf Y}}
\def\PROP{{\sc prop}} \def\hPROP{$\frac12${\sc prop}}
\def\bfk{{\bf k}}
\def\catPROP{{\tt PROP}} \def\cathPROP{{\textstyle\frac12}{\tt PROP}}
\def\catdgPROP{{\tt dgPROP}}\def\catho-dgPROP{\mbox{\tt ho-dgPROP}}
\def\Up{{\sf Up}}
\def\Dw{{\sf Dw}}
\def\span{{\it Span}}
\def\sigmaspan{{\it Span}_{\mbox{\scriptsize $\Sigma$-$\Sigma$}}}
\def\Span{\span}
\def\pa{\partial}
\def\ot{\otimes}
\def\op{\oplus}
\def\id{1\!\!1}
\def\Rada#1#2#3{#1_{#2},\dots,#1_{#3}}
\def\sfP{{\sf P}}
\def\sfD{{\sf D}}
\def\sfS{{\sf S}} \def\S{{\sfS}}

\def\doubless#1#2{{
\def\arraystretch{.5}
\begin{array}{c}
\mbox{\scriptsize $\scriptstyle #1$}
\\
\mbox{\scriptsize $\scriptstyle #2$}
\end{array}\def\arraystretch{1}
}}
\def\cases#1#2#3#4{
                  \left\{
                         \begin{array}{ll}
                           #1,\ &\mbox{#2}
                           \\
                           #3,\ &\mbox{#4}
                          \end{array}
                   \right.
}

\newcommand{\Vtriangle}[6]{
\setlength{\unitlength}{1.3em}
\begin{picture}(5,3.6)(0,-.1)
\thinlines
\put(0,2.5){\makebox(0,0){$#1$}}
\put(5,2.5){\makebox(0,0){$#2$}}
\put(2.5,0){\makebox(0,0){$#3$}}

\put(1,2.5){\vector(1,0){3}}
\put(0.5,2){\vector(1,-1){1.5}}
\put(4.5,2){\vector(-1,-1){1.5}}

\put(1,1){\makebox(0,0)[r]{\scriptsize $#5$}}
\put(2.5,3){\makebox(0,0)[b]{\scriptsize $#4$}}
\put(4,1){\makebox(0,0)[l]{\scriptsize $#6$}}
\end{picture}
}

\def\krouzek#1{
\setlength{\unitlength}{#1cm}
\bezier{100}(-2,0)(-2,.5)(-1.73,1)
\bezier{100}(-1.73,1)(-1.48,1.48)(-1,1.73)
\bezier{100}(-1,1.73)(-.5,2)(0,2)
\bezier{100}(2,0)(2,.5)(1.73,1)
\bezier{100}(1.73,1)(1.48,1.48)(1,1.73)
\bezier{100}(1,1.73)(.5,2)(0,2)
\bezier{100}(2,0)(2,-.5)(1.73,-1)
\bezier{100}(1.73,-1)(1.48,-1.48)(1,-1.73)
\bezier{100}(1,-1.73)(.5,-2)(0,-2)
\bezier{100}(-2,0)(-2,-.5)(-1.73,-1)
\bezier{100}(-1.73,-1)(-1.48,-1.48)(-1,-1.73)
\bezier{100}(-1,-1.73)(-.5,-2)(0,-2)
}

\def\dvojiteypsilon{{
\unitlength=.3pt
\begin{picture}(24.00,30.00)(0.00,3.00)
\put(10.00,20.00){\line(0,-1){10.00}}
\bezier{20}(10.00,10.00)(15,5)(20.00,0.00)
\bezier{20}(10.00,10.00)(5,5)(0.00,0.00)
\bezier{20}(10.00,20.00)(15,25)(20.00,30.00)
\bezier{20}(0.00,30.00)(5,25)(10.00,20.00)
\end{picture}}}

\def\motylek{{
\unitlength=.3pt
\begin{picture}(66.00,60.00)(-3.00,20.00)
\put(10.00,50.00){\line(0,1){0.00}}
\put(50.00,10.00){\line(0,-1){10.00}}
\put(10.00,10.00){\line(0,-1){10.00}}
\put(50.00,60.00){\line(0,-1){10.00}}
\put(10.00,60.00){\line(0,-1){10.00}}
\put(60.00,40.00){\line(0,-1){20.00}}
\put(0.00,40.00){\line(0,-1){20.00}}
\bezier{20}(50.00,10.00)(55.00,15.00)(60.00,20.00)
\bezier{20}(0.00,20.00)(5.00,15.00)(10.00,10.00)
\bezier{20}(50.00,50.00)(55.00,45.00)(60.00,40.00)
\bezier{20}(10.00,10.00)(15,15)(22,22)
\put(28,28){\bezier{20}(10.00,10.00)(15,15)(22,22)}
\put(10.00,50.00){\line(1,-1){40.00}}
\bezier{20}(10.00,50.00)(5.00,45.00)(0.00,40.00)
\end{picture}}}

\def\motylekmeziradky{{
\unitlength=.3pt
\begin{picture}(66.00,30.00)(-3.00,20.00)
\put(10.00,50.00){\line(0,1){0.00}}
\put(50.00,10.00){\line(0,-1){10.00}}
\put(10.00,10.00){\line(0,-1){10.00}}
\put(50.00,60.00){\line(0,-1){10.00}}
\put(10.00,60.00){\line(0,-1){10.00}}
\put(60.00,40.00){\line(0,-1){20.00}}
\put(0.00,40.00){\line(0,-1){20.00}}
\bezier{20}(50.00,10.00)(55.00,15.00)(60.00,20.00)
\bezier{20}(0.00,20.00)(5.00,15.00)(10.00,10.00)
\bezier{20}(50.00,50.00)(55.00,45.00)(60.00,40.00)
\bezier{20}(10.00,10.00)(15,15)(22,22)
\put(28,28){\bezier{20}(10.00,10.00)(15,15)(22,22)}
\put(10.00,50.00){\line(1,-1){40.00}}
\bezier{20}(10.00,50.00)(5.00,45.00)(0.00,40.00)
\end{picture}}}

\def\motylekmensi{{
\unitlength=.19pt
\begin{picture}(66.00,60.00)(-3.00,9.00)
\put(10.00,50.00){\line(0,1){0.00}}
\put(50.00,10.00){\line(0,-1){10.00}}
\put(10.00,10.00){\line(0,-1){10.00}}
\put(50.00,60.00){\line(0,-1){10.00}}
\put(10.00,60.00){\line(0,-1){10.00}}
\put(60.00,40.00){\line(0,-1){20.00}}
\put(0.00,40.00){\line(0,-1){20.00}}
\bezier{20}(50.00,10.00)(55.00,15.00)(60.00,20.00)
\bezier{20}(0.00,20.00)(5.00,15.00)(10.00,10.00)
\bezier{20}(50.00,50.00)(55.00,45.00)(60.00,40.00)
\bezier{20}(10.00,10.00)(15,15)(22,22)
\put(28,28){\bezier{20}(10.00,10.00)(15,15)(22,22)}
\bezier{40}(10.00,50.00)(30,30)(50,10)
\bezier{20}(10.00,50.00)(5.00,45.00)(0.00,40.00)
\end{picture}}}

\def\zeroatwo#1{%
{
\unitlength=.5pt
\begin{picture}(36.00,30.00)(-3.00,10.00)
\put(10.00,20.00){\makebox(0.00,0.00)[l]{$#1$}}
\put(20.00,10.00){\line(0,-1){5.00}}
\put(10.00,10.00){\line(0,-1){5.00}}
\put(30.00,30.00){\line(-1,0){30.00}}
\put(30.00,10.00){\line(0,1){20.00}}
\put(0.00,10.00){\line(1,0){30.00}}
\put(0.00,30.00){\line(0,-1){20.00}}
\end{picture}}
}

\def\twoazero#1{%
{
\unitlength=.5pt
\begin{picture}(36.00,30.00)(-3.00,0.00)
\put(20.00,25.00){\line(0,-1){5.00}}
\put(10.00,25.00){\line(0,-1){5.00}}
\put(10.00,10.00){\makebox(0.00,0.00)[l]{$#1$}}
\put(30.00,20.00){\line(-1,0){30.00}}
\put(30.00,0.00){\line(0,1){20.00}}
\put(0.00,0.00){\line(1,0){30.00}}
\put(0.00,20.00){\line(0,-1){20.00}}
\end{picture}}
}

\def\levaplastev{
\unitlength=.35pt
\begin{picture}(38,30)(-34,4)
\put(0.00,0.00){\line(0,1){10}}
\put(-20.00,0.00){\line(0,1){10}}
\put(-10.00,20.00){\line(0,1){10}}
\put(-30.00,20.00){\line(0,1){10}}
\put(0,10){\bezier{20}(0.00,0.00)(-5.00,5.00)(-10.00,10.00)}
\put(-20,10){\bezier{20}(0.00,0.00)(-5.00,5.00)(-10.00,10.00)}
\put(-20,10){\bezier{20}(0.00,0.00)(5.00,5.00)(10.00,10.00)}
\end{picture}
}

\def\twoaone#1{%
{
\unitlength=.5pt
\begin{picture}(36.00,30.00)(-3.00,-10.00)
\put(20.00,25.00){\line(0,-1){5.00}}
\put(15.00,0.00){\line(0,-1){5.00}}
\put(10.00,25.00){\line(0,-1){5.00}}
\put(10.00,10.00){\makebox(0.00,0.00)[l]{$#1$}}
\put(30.00,20.00){\line(-1,0){30.00}}
\put(30.00,0.00){\line(0,1){20.00}}
\put(0.00,0.00){\line(1,0){30.00}}
\put(0.00,20.00){\line(0,-1){20.00}}
\end{picture}}
}
\def\zeroathree#1{
{
\unitlength=.5pt
\begin{picture}(48.00,30.00)(-4.00,10.00)
\put(20.00,20.00){\makebox(0.00,0.00){$#1$}}
\put(0.00,30.00){\line(0,-1){20.00}}
\put(40.00,30.00){\line(-1,0){40.00}}
\put(40.00,10.00){\line(0,1){20.00}}
\put(0.00,10.00){\line(1,0){40.00}}
\put(30.00,10.00){\line(0,-1){5.00}}
\put(20.00,10.00){\line(0,-1){5.00}}
\put(10.00,10.00){\line(0,-1){5.00}}
\end{picture}}
}

\def\twoatwo#1{
{
\unitlength=.5pt
\begin{picture}(30.00,40.00)(0.00,0.00)
\put(10.00,20.00){\makebox(0.00,0.00)[l]{$#1$}}
\put(20.00,10.00){\line(0,-1){5.00}}
\put(10.00,10.00){\line(0,-1){5.00}}
\put(20.00,35.00){\line(0,-1){5.00}}
\put(10.00,35.00){\line(0,-1){5.00}}
\put(0.00,10.00){\line(0,1){20.00}}
\put(30.00,10.00){\line(-1,0){30.00}}
\put(30.00,30.00){\line(0,-1){20.00}}
\put(0.00,30.00){\line(1,0){30.00}}
\end{picture}}
}

\def\jednadva{{
\unitlength=.4pt
\begin{picture}(24.00,20.00)(-2.00,0.00)
\bezier{20}(10.00,10.00)(15.00,5.00)(20.00,0.00)
\bezier{20}(10.00,10.00)(5.00,5.00)(0.00,0.00)
\put(10.00,20.00){\line(0,-1){10.00}}
\end{picture}}
}

\def\jednactyri{{
\unitlength=.05pt
\begin{picture}(176.00,160.00)(-8.00,0.00)
\put(80.00,100.00){\line(0,1){60.00}}
\bezier{20}(80.00,80.00)(100.00,30.00)(110.00,0.00)
\bezier{20}(80.00,80.00)(60.00,30.00)(50.00,0.00)
\bezier{20}(80.00,80.00)(120.00,40.00)(160.00,0.00)
\bezier{20}(80.00,80.00)(40.00,40.00)(0.00,0.00)
\put(80.00,100.00){\line(0,-1){20.00}}
\end{picture}}
}

\def\ctyrijedna{{
\unitlength=.05pt
\begin{picture}(176.00,160.00)(-8.00,-160.00)
\put(80.00,-100.00){\line(0,-1){60.00}}
\bezier{20}(80.00,-80.00)(100.00,-30.00)(110.00,0.00)
\bezier{20}(80.00,-80.00)(60.00,-30.00)(50.00,0.00)
\bezier{20}(80.00,-80.00)(120.00,-40.00)(160.00,0.00)
\bezier{20}(80.00,-80.00)(40.00,-40.00)(0.00,0.00)
\put(80.00,-80.00){\line(0,-1){40.00}}
\end{picture}}
}

\def\dvajedna{{
\unitlength=.4pt
\begin{picture}(24.00,20.00)(-2.00,0.00)
\put(10.00,10.00){\line(0,-1){10.00}}
\bezier{20}(10.00,10.00)(15.00,15.00)(20.00,20.00)
\bezier{20}(0.00,20.00)(5.00,15.00)(10.00,10.00)
\end{picture}}
}

\def\dvadva{{
\unitlength=.8pt
\begin{picture}(12.00,10.00)(-1.00,0.00)
\bezier{30}(0.00,0.00)(5.00,5.00)(10.00,10.00)
\bezier{30}(0.00,10.00)(5.00,5.00)(10.00,0.00)
\end{picture}}
}

\def\jednatri{{
\unitlength=.4pt
\begin{picture}(24.00,20.00)(-2.00,0.00)
\bezier{20}(10.00,10.00)(15.00,5.00)(20.00,0.00)
\bezier{20}(10.00,10.00)(5.00,5.00)(0.00,0.00)
\put(10.00,20.00){\line(0,-1){20.00}}
\end{picture}}
}

\def\trijedna{{
\unitlength=.4pt
\begin{picture}(24.00,20.00)(-2.00,-20.00)
\bezier{20}(10.00,-10.00)(15.00,-5.00)(20.00,0.00)
\bezier{20}(10.00,-10.00)(5.00,-5.00)(0.00,0.00)
\put(10.00,-20.00){\line(0,1){20.00}}
\end{picture}}
}

\def\dvatri{{
\unitlength=0.4pt
\begin{picture}(24.00,20.00)(-2.00,0.00)
\put(10.00,10.00){\line(0,-1){10.00}}
\bezier{30}(0.00,0.00)(10.00,10.00)(20.00,20.00)
\bezier{30}(0.00,20.00)(10.00,10.00)(20.00,0.00)
\end{picture}}
}

\def\tridva{{
\unitlength=.4pt
\begin{picture}(24.00,20.00)(-2.00,-20.00)
\put(10.00,-10.00){\line(0,1){10.00}}
\bezier{30}(0.00,0.00)(10.00,-10.00)(20.00,-20.00)
\bezier{30}(0.00,-20.00)(10.00,-10.00)(20.00,0.00)
\end{picture}}
}

\def\tritri{{
\unitlength=.4pt
\begin{picture}(24.00,20.00)(-2.00,0.00)
\bezier{30}(0.00,0.00)(10.00,10.00)(20.00,20.00)
\bezier{30}(0.00,20.00)(10.00,10.00)(20.00,0.00)
\put(10.00,20.00){\line(0,-1){20.00}}
\end{picture}}
}

\def\dvactyri{{
\unitlength=.1pt
\begin{picture}(96.00,80.00)(-8.00,0.00)
\bezier{30}(0.00,0.00)(40.00,40.00)(80.00,80.00)
\bezier{30}(0.00,80.00)(40.00,40.00)(80.00,0.00)
\bezier{20}(40.00,40.00)(48.50,14.00)(53.50,0.00)
\bezier{20}(40.00,40.00)(32.50,16.00)(26.00,0.00)
\end{picture}}
}

\def\dvaZbbbZb{{
\unitlength=0.1pt
\begin{picture}(96.00,80.00)(-8.00,0.00)
\bezier{10}(20.00,20.00)(30.00,10.00)(40.00,0.00)
\put(20.00,20.00){\line(0,-1){20.00}}
\bezier{30}(0.00,0.00)(40.00,40.00)(80.00,80.00)
\bezier{30}(0.00,80.00)(40.00,40.00)(80.00,0.00)
\end{picture}}
}

\def\dvabZbbbZ{{
\unitlength= 0.1pt
\begin{picture}(96.00,80.00)(-8.00,0.00)
\bezier{10}(60.00,20.00)(50.00,10.00)(40.00,0.00)
\put(60.00,20.00){\line(0,-1){20.00}}
\bezier{30}(0.00,0.00)(40.00,40.00)(80.00,80.00)
\bezier{30}(0.00,80.00)(40.00,40.00)(80.00,0.00)
\end{picture}}
}

\def\dvaZbbZbb{{
\unitlength=0.1pt
\begin{picture}(96.00,80.00)(-8.00,0.00)
\bezier{8}(11.00,9.00)(15.50,4.00)(21.00,0.00)
\put(40.00,40.00){\line(0,-1){40.00}}
\bezier{30}(0.00,0.00)(40.00,40.00)(80.00,80.00)
\bezier{30}(0.00,80.00)(40.00,40.00)(80.00,0.00)
\end{picture}}
}

\def\dvabbZbbZ{{
\unitlength=0.1pt
\begin{picture}(96.00,80.00)(-8.00,0.00)
\bezier{5}(70.00,10.00)(65.00,5.50)(60.00,0.00)
\put(40.00,40.00){\line(0,-1){40.00}}
\bezier{30}(0.00,0.00)(40.00,40.00)(80.00,80.00)
\bezier{30}(0.00,80.00)(40.00,40.00)(80.00,0.00)
\end{picture}}
}

\def\dvabZbbZb{{
\unitlength=0.1pt
\begin{picture}(96.00,80.00)(-8.00,0.00)
\bezier{10}(40.00,20.00)(50.00,10.00)(60.00,0.00)
\bezier{10}(40.00,20.00)(30.00,10.00)(20.00,0.00)
\put(40.00,40.00){\line(0,-1){20.00}}
\bezier{30}(0.00,0.00)(40.00,40.00)(80.00,80.00)
\bezier{30}(0.00,80.00)(40.00,40.00)(80.00,0.00)
\end{picture}}
}

\def\ZbbZbZb{{
\unitlength=.05pt
\begin{picture}(192.00,120.00)(-16.00,0.00)
\bezier{10}(30.00,30.00)(50.00,10.00)(60.00,0.00)
\bezier{20}(50.00,50.00)(70.00,30.00)(100.00,0.00)
\bezier{30}(80.00,80.00)(120.00,40.00)(160.00,0.00)
\bezier{30}(80.00,80.00)(40.00,40.00)(0.00,0.00)
\put(80.00,160.00){\line(0,-1){80.00}}
\end{picture}}
}

\def\dvacarkatri{{
\unitlength=.2pt
\begin{picture}(40.00,50.00)(0.00,0.00)
\put(20.00,30.00){\line(0,-1){10.00}}
\put(20.00,0.00){\line(0,1){20}}
\bezier{20}(20.00,20.00)(30.00,10.00)(40.00,0.00)
\bezier{20}(20.00,20.00)(10.00,10.00)(0.00,0.00)
\bezier{20}(20.00,30.00)(30.00,40.00)(40.00,50.00)
\bezier{20}(0.00,50.00)(10.00,40.00)(20.00,30.00)
\end{picture}}
}

\def\dvacarkaZbbZb{{
\unitlength=.2pt
\begin{picture}(40.00,50.00)(0.00,0.00)
\put(20.00,30.00){\line(0,-1){10.00}}
\bezier{10}(20.00,0.00)(15.00,5.00)(10.00,10.00)
\bezier{20}(20.00,20.00)(30.00,10.00)(40.00,0.00)
\bezier{20}(20.00,20.00)(10.00,10.00)(0.00,0.00)
\bezier{20}(20.00,30.00)(30.00,40.00)(40.00,50.00)
\bezier{20}(0.00,50.00)(10.00,40.00)(20.00,30.00)
\end{picture}}
}

\def\dvacarkabZbbZ{{
\unitlength=.2pt
\begin{picture}(40.00,50.00)(0.00,0.00)
\put(20.00,30.00){\line(0,-1){10.00}}
\bezier{10}(20.00,0.00)(25.00,5.00)(30.00,10.00)
\bezier{20}(20.00,20.00)(30.00,10.00)(40.00,0.00)
\bezier{20}(20.00,20.00)(10.00,10.00)(0.00,0.00)
\bezier{20}(20.00,30.00)(30.00,40.00)(40.00,50.00)
\bezier{20}(0.00,50.00)(10.00,40.00)(20.00,30.00)
\end{picture}}
}

\def\bZbZbbZ{{
\unitlength=0.05pt
\begin{picture}(192.00,160.00)(-16.00,0.00)
\bezier{10}(100.00,0.00)(120.00,20.00)(130.00,30.00)
\bezier{20}(110.00,50.00)(90.00,30.00)(60.00,0.00)
\put(80.00,120.00){\line(0,1){40.00}}
\bezier{30}(80.00,80.00)(120.00,40.00)(160.00,0.00)
\bezier{30}(80.00,80.00)(40.00,40.00)(0.00,0.00)
\put(80.00,120.00){\line(0,-1){40.00}}
\end{picture}}
}

\def\gen#1#2{
\if #11
    \if #22 \jednadva \else \fi
\else
\fi
\if #12
    \if #22 \dvadva \else \fi
\else
\fi
\if #12
    \if #21 \dvajedna \else \fi
\else
\fi
\if #13
    \if #22 \tridva \else \fi
\else
\fi
\if #13
    \if #21 \trijedna \else \fi
\else
\fi
\if #12
    \if #23 \dvatri \else \fi
\else
\fi
\if #11
    \if #23 \jednatri \else \fi
\else
\fi
\if #11
    \if #24 \jednactyri \else \fi
\fi
\if #14
    \if #21 \ctyrijedna \else \fi
\fi
}

\def\onethree{{
\unitlength=0.3pt
\begin{picture}(44.00,40.00)(-2.00,0.00)
\put(30.00,10.00){\line(0,-1){10.00}}
\put(20.00,10.00){\line(0,-1){10.00}}
\put(10.00,10.00){\line(0,-1){10.00}}
\put(20.00,40.00){\line(0,-1){10.00}}
\put(0.00,10.00){\line(0,1){20.00}}
\put(40.00,10.00){\line(-1,0){40.00}}
\put(40.00,30.00){\line(0,-1){20.00}}
\put(0.00,30.00){\line(1,0){40.00}}
\end{picture}}
}

\def\twotwo{{
\unitlength=.18pt
\begin{picture}(44.00,50.00)(-2.00,0.00)
\put(0.00,10.00){\line(0,1){30.00}}
\put(40.00,10.00){\line(-1,0){40.00}}
\put(40.00,40.00){\line(0,-1){30.00}}
\put(0.00,40.00){\line(1,0){40.00}}
\put(30.00,10.00){\line(0,-1){10.00}}
\put(10.00,10.00){\line(0,-1){10.00}}
\put(30.00,50.00){\line(0,-1){10.00}}
\put(10.00,50.00){\line(0,-1){10.00}}
\end{picture}}
}

\def\bZbbZ{
{
\unitlength=.27pt
\begin{picture}(48.00,30.00)(-4,0.00)
\bezier{34}(20.00,20.00)(30.00,10.00)(40.00,0.00)
\bezier{34}(20.00,20.00)(10.00,10.00)(0.00,0.00)
\bezier{20}(30.00,10.00)(25.00,5.00)(20.00,0.00)
\put(20.00,30.00){\line(0,-1){10.00}}
\end{picture}} 
}

\def\ZbbZb{{
\unitlength=.27pt
\begin{picture}(48.00,30.00)(-4,0.00)
\bezier{34}(20.00,20.00)(30.00,10.00)(40.00,0.00)
\bezier{34}(20.00,20.00)(10.00,10.00)(0.00,0.00)
\bezier{20}(10.00,10.00)(15.00,5.00)(20.00,0.00)
\put(20.00,30.00){\line(0,-1){10.00}}
\end{picture}}}

\def\ZbbZbb{{
\unitlength=.07pt
\begin{picture}(192.00,120.00)(-16.00,0.00)
\bezier{20}(20.00,20.00)(30.00,10.00)(40.00,0.00)
\bezier{30}(80.00,80.00)(120.00,40.00)(160.00,0.00)
\bezier{30}(80.00,80.00)(40.00,40.00)(0.00,0.00)
\put(80.00,120.00){\line(0,-1){120.00}}
\end{picture}}
}

\def\bbZbbZ{{
\unitlength=0.07pt
\begin{picture}(192.00,120.00)(-16.00,0.00)
\bezier{20}(140.00,20.00)(130.00,10.00)(120.00,0.00)
\bezier{30}(80.00,80.00)(120.00,40.00)(160.00,0.00)
\bezier{30}(80.00,80.00)(40.00,40.00)(0.00,0.00)
\put(80.00,120.00){\line(0,-1){120.00}}
\end{picture}}
}

\def\bZbbZb{{
\unitlength=.07pt
\begin{picture}(192.00,120.00)(-16.00,0.00)
\bezier{20}(80.00,20.00)(90.00,10.00)(100.00,0.00)
\bezier{20}(80.00,20.00)(70.00,10.00)(60.00,0.00)
\put(80.00,40.00){\line(0,-1){20.00}}
\put(80.00,40.00){\line(0,1){0.00}}
\put(80.00,80.00){\line(0,-1){40.00}}
\put(80.00,120.00){\line(0,-1){40.00}}
\bezier{30}(80.00,80.00)(120.00,40.00)(160.00,0.00)
\bezier{30}(80.00,80.00)(40.00,40.00)(0.00,0.00)
\end{picture}}
}

\def\ZbbbZb{{
\unitlength=.07pt
\begin{picture}(192.00,120.00)(-16.00,0.00)
\put(40.00,0.00){\line(0,1){20.00}}
\put(40.00,40.00){\line(0,-1){40.00}}
\bezier{20}(40.00,40.00)(60.00,20.00)(80.00,0.00)
\put(80.00,40.00){\line(0,1){0.00}}
\put(80.00,120.00){\line(0,-1){40.00}}
\bezier{30}(80.00,80.00)(120.00,40.00)(160.00,0.00)
\bezier{30}(80.00,80.00)(40.00,40.00)(0.00,0.00)
\end{picture}}
}

\def\bZbbbZ{{
\unitlength=0.07pt
\begin{picture}(192.00,120.00)(-16.00,0.00)
\put(120.00,40.00){\line(0,-1){40.00}}
\bezier{20}(120.00,40.00)(100.00,20.00)(80.00,0.00)
\put(80.00,40.00){\line(0,1){0.00}}
\put(80.00,120.00){\line(0,-1){40.00}}
\bezier{30}(80.00,80.00)(120.00,40.00)(160.00,0.00)
\bezier{30}(80.00,80.00)(40.00,40.00)(0.00,0.00)
\end{picture}}
}

\def\ZvvZv{{
\unitlength=.27pt
\begin{picture}(48.00,30.00)(-4,-30.00)
\bezier{30}(20.00,-20.00)(30.00,-10.00)(40.00,0.00)
\bezier{30}(20.00,-20.00)(10.00,-10.00)(0.00,0.00)
\bezier{20}(10.00,-10.00)(15.00,-5.00)(20.00,0.00)
\put(20.00,-30.00){\line(0,1){10.00}}
\end{picture}}}

\def\vZvvZ{{
\unitlength=.27pt
\begin{picture}(48.00,30.00)(-4,-30.00)
\bezier{34}(20.00,-20.00)(30.00,-10.00)(40.00,0.00)
\bezier{34}(20.00,-20.00)(10.00,-10.00)(0.00,0.00)
\bezier{20}(30.00,-10.00)(25.00,-5.00)(20.00,0.00)
\put(20.00,-30.00){\line(0,1){10.00}}
\end{picture}}
}

\def\vZvvZdva{{
\unitlength=.2pt
\begin{picture}(48.00,40.00)(-4.00,0.00)
\bezier{10}(30.00,30.00)(25.00,35.00)(20.00,40.00)
\bezier{34}(0.00,0.00)(20.00,20.00)(40.00,40.00)
\bezier{34}(0.00,40.00)(20.00,20.00)(40.00,0.00)
\end{picture}}}

\def\dvabZbbZ{{
\unitlength=.2pt
\begin{picture}(48.00,40.00)(-4.00,-40.00)
\bezier{10}(30.00,-30.00)(25.00,-35.00)(20.00,-40.00)
\bezier{34}(0.00,0.00)(20.00,-20.00)(40.00,-40.00)
\bezier{34}(0.00,-40.00)(20.00,-20.00)(40.00,0.00)
\end{picture}}}

\def\triZbbZb{{
\unitlength=.2pt
\begin{picture}(48.00,40.00)(-4.00,0.00)
\put(20.00,40.00){\line(0,-1){20.00}}
\bezier{10}(10.00,10.00)(15.00,5.00)(20.00,0.00)
\bezier{34}(0.00,0.00)(20.00,20.00)(40.00,40.00)
\bezier{34}(0.00,40.00)(20.00,20.00)(40.00,0.00)
\end{picture}}}

\def\ZvvZvdva{{
\unitlength=.2pt
\begin{picture}(48.00,40.00)(-4.00,0.00)
\bezier{10}(10.00,30.00)(15.00,35.00)(20.00,40.00)
\bezier{34}(0.00,0.00)(20.00,20.00)(40.00,40.00)
\bezier{34}(0.00,40.00)(20.00,20.00)(40.00,0.00)
\end{picture}}}

\def\tribZbbZ{{
\unitlength=.2pt
\begin{picture}(48.00,40.00)(-4.00,0.00)
\put(20.00,40.00){\line(0,-1){20.00}}
\bezier{10}(30.00,10.00)(25.00,5.00)(20.00,0.00)
\bezier{34}(0.00,0.00)(20.00,20.00)(40.00,40.00)
\bezier{34}(0.00,40.00)(20.00,20.00)(40.00,0.00)
\end{picture}}
}

\def\ZvvZvtri{{
\unitlength=.2pt
\begin{picture}(48.00,40.00)(-4.00,0.00)
\put(20.00,0.00){\line(0,1){20.00}}
\bezier{10}(10.00,30.00)(15.00,35.00)(20.00,40.00)
\bezier{34}(0.00,0.00)(20.00,20.00)(40.00,40.00)
\bezier{34}(0.00,40.00)(20.00,20.00)(40.00,0.00)
\end{picture}}
}

\def\vZvvZtri{{
\unitlength=.2pt
\begin{picture}(48.00,40.00)(-4.00,0.00)
\put(20.00,0.00){\line(0,1){20.00}}
\bezier{10}(30.00,30.00)(25.00,35.00)(20.00,40.00)
\bezier{34}(0.00,0.00)(20.00,20.00)(40.00,40.00)
\bezier{34}(0.00,40.00)(20.00,20.00)(40.00,0.00)
\end{picture}}
}

\def\dvaZbbZb{{
\unitlength=.2pt
\begin{picture}(48.00,40.00)(-4.00,-40.00)
\bezier{10}(10.00,-30.00)(15.00,-35.00)(20.00,-40.00)
\bezier{34}(0.00,0.00)(20.00,-20.00)(40.00,-40.00)
\bezier{34}(0.00,-40.00)(20.00,-20.00)(40.00,0.00)
\end{picture}}}

\title[Minimal model of the PROP for bialgebras]%
      {A resolution (minimal model) of the PROP for bialgebras}

\author[M. Markl]{Martin Markl}
\thanks{The author was supported by the
grant GA AV \v CR \#1019203. Preliminary results were announced at 
Workshop on Topology, Operads \&
Quantization, Warwick, UK, 11.12.~2001.}

\begin{abstract}
This paper is concerned with a minimal resolution of the \PROP\ for
bialgebras (Hopf algebras without unit, counit and
antipode). We prove a theorem about the form of this resolution
(Theorem~\ref{main}) and give, in Section~\ref{jeden_medvidek}, 
a lot of explicit formulas for the differential.
\end{abstract}

\def\subjclassname{\textup{2000} Mathematics Subject
     Classification}

\subjclass[2000]{16W30, 57T05, 18C10, 18G99}
\email{markl@math.cas.cz}
\address{Mathematical Institute of the Academy, {\v Z}itn{\'a} 25, 
         115 67 Prague 1, The Czech Republic}
\date{January 5, 2005}

\maketitle

\section{Introduction and main results}
\label{intro}

A {\em bialgebra\/} is a vector space $V$ with 
a {\em multiplication\/} $\mu : V \ot V \to V$ and a {\em
comultiplication\/} (also called a {\em diagonal\/}) $\Delta : V \to
V\ot V$. The multiplication is associative:
\begin{equation}
\label{ass}
\mu(\mu \ot \id_V) = \mu(\id_V \ot \mu), 
\end{equation}
where $\id_V : V \to V$ denotes the identity map,
the comultiplication is coassociative:
\begin{equation}
\label{coass}
(\id_V \ot \Delta)\Delta = (\Delta \ot \id_V)\Delta
\end{equation}
and the usual compatibility relation between $\mu$ and $\Delta$ is assumed:
\begin{equation}
\label{compatibility}
\Delta \circ \mu = (\mu \ot \mu) T_{\sigma(2,2)} (\Delta \ot \Delta),
\end{equation}
where $T_{\sigma(2,2)}: V^{\ot 4} \to V^{\ot 4}$ is defined by  
\[
T_{\sigma(2,2)}(v_1 \ot v_2 \ot v_3 \ot v_4) := v_1 \ot v_3 \ot v_2 \ot v_4,
\]
for $v_1,v_2,v_3,v_4 \in V$ (the meaning of the notation $\sigma(2,2)$
will be explained in Definition~\ref{1}). 
We suppose that $V$, as well as all other
algebraic objects in this paper, are defined over a field $\bfk$ of
characteristic zero.

Let $\sfB$ be the $\bfk$-linear \PROP\ (see~\cite{markl:ws93,markl:JPAA96} 
or Section~\ref{Furby} of this paper for
the terminology) describing bialgebras. The goal of this 
paper is to describe a {\em minimal model\/} of
$\sfB$, that is, a differential graded (dg) 
{\bfk}-linear {\sc prop} $({\sf M},\pa)$
together with a homology isomorphism
\[
(\sfB,0) \stackrel{\rho}{\longleftarrow} ({\sf M},\pa)
\] 
such that 
\begin{itemize}
\item[(i)]
the \PROP\ ${\sf M}$ is free and 
\item[(ii)]
the image of $\pa$ consists of decomposable elements of ${\sf M}$ (the
minimality condition),
\end{itemize}
see again Section~\ref{Furby} where free \PROP{s} and decomposable
elements are recalled.

The initial stages of this minimal model were constructed
in~\cite[page~145]{markl:ws93} and
\cite[pages~215--216]{markl:JPAA96}. According to our general
philosophy, it should contain all information about the deformation
theory of bialgebras.  In particular, the Gerstenhaber-Schack
cohomology which is known to control deformations of
bialgebras~\cite{gerstenhaber-schack:Proc.Nat.Acad.Sci.USA90} can be
read off from this model as follows.

Let $\End_V$ denote the endomorphism \PROP\ of $V$ and let a bialgebra
structure $B = (V,\mu,\Delta)$ on $V$ be given by a homomorphism of
\PROP{s} $\beta : \sfB \to \End_V$. The composition $\beta \circ \rho
: \sfM \to \End_V$ makes $\End_V$ an $\sfM$-module 
(in the sense of~\cite[page~203]{markl:JPAA96}), 
therefore one may consider the vector space of derivations $\DER$. For $\theta
\in \DER$ define $\delta \theta : = \theta \circ \pa$. It follows from
the obvious fact that $\rho \circ \pa = 0$ that $\delta \theta$ is again
a derivation, so $\delta$ is a well-defined endomorphism of the
vector space $\DER$ which clearly satisfies $\delta^2=0$. Then
\[
H_b(B;B) \cong H(\DER,\delta),
\]
where $H_b(B;B)$ denotes the Gerstenhaber-Schack cohomology of the
bialgebra $B$ with coefficients in itself.

Algebras (in the sense recalled in Section~\ref{Furby}) over $({\sf
M},\pa)$ have all rights to be called {\em strongly homotopy
bialgebras\/}, that is, homotopy invariant versions of bialgebras, as
follows from principles explained in the introduction
of~\cite{markl:ha}.  This would mean, among other things, that, given
a structure of a dg-bialgebra on a chain complex $C_*$, then any chain
complex $D_*$, chain homotopy equivalent to $C_*$, has, in a certain
sense, a natural and unique structure of an algebra over our minimal
model $({\sf M},\pa)$.

For a discussion of \PROP{s} for bialgebras from another perspective,
see~\cite{pirashvili}. Constructions of various other (non-minimal)
resolutions of the \PROP\ for bialgebras, based mostly on a dg-version
of the Boardman-Vogt $W$-construction, will be the subject
of~\cite{kontsevich-soibelman:hopf-algebras}. A completely different
approach to bialgebras and resolutions of objects governing them can
be found in a series of papers by
Shoikhet~\cite{shoikhet:2/3,shoikhet,shoiket:explicit}, and also in a
recent draft by Saneblidze and Umble~\cite{saneblidze-umble:matrons}.
A general theory of resolutions of \PROP{s} is, besides~\cite{mv},
also the subject of Vallette's thesis and its
follow-up~\cite{vallette,vallette:CMR04}.

Let us briefly sketch the strategy of the construction of our model.
Consider objects $(V,\mu,\Delta)$, where $\mu: V \ot V \to V$ is 
an associative multiplication as in~(\ref{ass}), 
$\Delta: V \to V \ot V$ is a coassociative comultiplication as
in~(\ref{coass}), but the compatibility relation~(\ref{compatibility}) is
replaced by 
\begin{equation}
\label{pf}
\Delta \circ \mu = 0.
\end{equation}

\begin{definition}
A {\em half-bialgebra\/} or briefly {\em $\frac12$bialgebra\/} is a
vector space $V$ equipped with a multiplication $\mu$ and a
comultiplication $\Delta$ satisfying~(\ref{ass}),~(\ref{coass})
and~(\ref{pf}).
\end{definition}

We chose this strange name because~(\ref{pf}) is indeed, in a sense,
one half of the compatibility relation~(\ref{compatibility}).
{}For a formal variable $\epsilon$, consider the axiom
\[
\Delta \circ \mu = 
\epsilon \cdot (\mu \ot \mu) T_{\sigma(2,2)} (\Delta \ot \Delta).
\]
At $\epsilon = 1$ we get the usual
compatibility relation~(\ref{compatibility}) between the multiplication and the
diagonal, while $\epsilon =0$ gives~(\ref{pf}). 
Therefore~(\ref{compatibility}) can be interpreted
as a perturbation of~(\ref{pf}) which may be informally expressed
by saying that bialgebras are perturbations of $\frac12$bialgebras.
Experience with homological perturbation theory~\cite{halperin-stasheff}
leads us to formulate:

\vskip 3mm
\noindent
{\bf Principle.}
{\em\ The {\sc prop} $\sfB$ for bialgebras is a
perturbation of the \PROP\ $\sfhB$ for $\frac12$bialgebras. Therefore
there exists
a minimal model of the \PROP{} $\sfB$ that is 
a perturbation of a minimal model of the \PROP\ $\sfhB$
for $\frac12$bialgebras.
}

We therefore need to know a minimal model for $\sfhB$.  In general,
\PROP{s} are extremely huge objects, difficult to work with, but
$\frac12$bialgebras exist over much smaller objects than
\PROP{}s. These smaller objects, which we call \hPROP{s}, were
introduced in an e-mail message from
M.~Kontsevich~\cite{kontsevich:message} who called them small
\PROP{s}. The concept of \hPROP{s} makes the construction of a minimal
model of $\sfhB$ easy. We thus proceed in two steps.

{\bf Step 1.} We construct a minimal
model $(\Gamma(\Xi),\pa_0)$ of the \PROP\ $\sfhB$ for
$\frac12$bialgebras. Here $\Gamma(\Xi)$ denotes the free \PROP\ 
on the space of generators $\Xi$, see Theorem~\ref{dnes_tatarak_u_Jany}.

{\bf Step 2.} Our minimal model $({\sf M},\pa)$ of the \PROP\ $\sfB$ for
bialgebras will be then a perturbation of $(\Gamma(\Xi),\pa_0)$, that is, 
\[
({\sf M},\pa)  = (\Gamma(\Xi),\pa_0 + \papert),
\]
see Theorem~\ref{main}.

\vskip 1em

{\bf Acknowledegment.} I would like to express my gratitude to Jim
Stasheff, Steve Shnider, Vladimir Hinich, Wee Liang Gan and Petr
Somberg for careful reading the manuscript and many useful
suggestions. I would also like to thank the Erwin Schr\"odinger
International Institute for Mathematical Physics, Vienna, for the
hospitality during the period when the first draft of this paper was
completed.

My particular thanks are due to M.~Kontsevich whose
e-mail~\cite{kontsevich:message} shed a new light on the present work
and stimulated a cooperation with A.A.~Voronov which resulted in~\cite{mv}.
Also the referee's remarks were extremely helpful.

\vfill
\break

\section{Structure of \PROP{s} and \hPROP{s}}
\label{Furby}\label{sec2}

Let us recall that a $\bfk$-linear \PROP{} $\sfA$ (called a {\em
theory\/} in~\cite{markl:ws93,markl:JPAA96}) is a sequence of
$\bfk$-vector spaces $\{\sfA(m,n)\}_{m,n \geq 1}$ with compatible left
$\Sigma_m$- right $\Sigma_n$-actions and two types of equivariant
compositions, vertical:
\[
\circ: 
\sfA(m,u) \ot_{\Sigma_u} \sfA(u,n) \to \sfA(m,n), \ m,n,u \geq 1,
\] 
and horizontal:
\[
\boxtimes : \sfA(m_1,n_1) \otimes  \sfA(m_2,n_2) \to \sfA(m_1+m_2,n_1
+n_2),\
m_1,m_2,n_1,n_2 \geq 1,
\]
together with an identity $\id \in \sfA(1,1)$.  \PROP{s} should
satisfy axioms which could be read off from the example of the {\em
endomorphism \PROP\/} $\End_V$ of a vector space $V$, with
$\End_V(m,n)$ the space of linear maps ${\it Hom}_\bfk(\otexp Vn,\otexp
Vm)$, $\id \in \End_V(1,1)$ the identity map, horizontal composition
given by the tensor product of linear maps, and vertical composition
by the ordinary composition of maps.  One can therefore imagine
elements of $\sfA(m,n)$ as `abstract' maps with $n$ inputs and $m$
outputs.  See~\cite{maclane:RiceUniv.Studies63,markl:JPAA96} for
precise definitions.

We say that $X$ has {\em biarity\/} $(m,n)$ if $X \in \sfA(m,n)$.
We will sometimes use the operadic notation: for $X \in \sfA(m,k)$,
$Y \in \sfA(1,l)$ and $1 \leq i \leq k$, we write
\begin{equation}
\label{Zbraslavice1}
X \circ_i Y := X \circ (\id^{\ot (i-1)} \ot Y \ot  \id^{\ot (k-i)})
\in \sfA(m,k+l-1)
\end{equation}
and, similarly, for $U \in \sfA(k,1)$, $V \in \sfA(l,n)$ and $1 \leq j
\leq l$ we denote
\begin{equation}
\label{Zbraslavice2}
U \cric j V :=  (\id^{\ot (j-1)} \ot U \ot  \id^{\ot (l-j)}) \circ V
\in \sfA(k+l-1,n). 
\end{equation}

In~\cite{markl:JPAA96} we called a sequence $E = \{E(m,n)\}_{m,n \geq
1}$ of left $\Sigma_m$-, right $\Sigma_n$-$\bfk$-bimodules a {\em
core\/}, but we prefer now to call such sequences {\em
$\Sigma$-bimodules\/}. For any such a $\Sigma$-bimodule $E$, there
exists the {\em free \PROP\/} $\Gamma(E)$ generated by $E$. It also
makes sense to speak, in the category of \PROP{s}, about ideals,
presentations, modules,~etc, see~\cite[Chapter~2]{vallette} for details.

Recall that an {\em algebra over a \PROP\/} $\sfA$ is (given by) a
\PROP\ morphism $\alpha : \sfA \to \End_V$.  A \PROP\ $\sfA$ is {\em
augmented\/} if there exist a homomorphism $\epsilon : \sfA \to {\sf
1}$ (the {\em augmentation\/}) to the trivial \PROP\ ${\sf 1} :=
\End_\bfk$.  Therefore an augmentation is the same as a structure of
an $\sfA$-algebra on the one-dimensional vector space $\bfk$.

Let $\sfA^+ := {\rm Ker}(\epsilon)$ denote the {\em augmentation
ideal} of an augmented \PROP\ $\sfA$. The space $D(\sfA) := \sfA^+
\hskip -.2mm \circ \sfA^+$ is then called the {\em space of
decomposables\/} and the quotient $Q(\sfA) := \sfA^+/D(\sfA)$ the {\em
space of indecomposables\/} of the augmented \PROP\ $\sfA$.
Observe that each free \PROP\ $\Gamma(E)$ is canonically augmented,
with the augmentation defined by $\epsilon(E) := 0$.

Let $\Gamma(\gen12,\gen21)$ be the free \PROP\ generated by one
operation $\gen12$ of biarity $(1,2)$ and one operation $\gen21$ of
biarity $(2,1)$. More formally, $\Gamma(\gen12,\gen21) := \Gamma(E)$
with $E$ the $\Sigma$-bimodule $\bfk \cdot \gen12 \ot \bfk[\Sigma_2]
\oplus \bfk[\Sigma_2] \otimes \bfk \cdot \gen21$.  As we explained
in~\cite{markl:ws93,markl:JPAA96}, the \PROP\ $\sfB$ describing
bialgebras has a presentation
\begin{equation}
\label{pekne_jsem_si_ve_Zbraslavicich_poletal}
\sfB = \Gamma(\gen12,\gen21)/{\sf I}_{\sf B},
\end{equation}
where ${\sf I}_{\sf B}$ denotes the ideal generated by
\[
\ZbbZb - \bZbbZ,\
\ZvvZv - \vZvvZ\  \mbox { and }\ 
\dvojiteypsilon - \motylek \hskip .2em.
\]
In the above display we denoted
\begin{eqnarray*}
&\ZbbZb := \gen12(\gen12 \ot \id),\
\bZbbZ := \gen12(\id \ot \gen12),\
\ZvvZv := (\gen21 \ot \id)\gen21,\
\vZvvZ := (\id \ot \gen21)\gen21,&
\\
&\dvojiteypsilon  := \gen21 \circ \gen12
\mbox {\hskip .5em and \hskip .6em}
\motylek := (\gen12 \ot \gen12)\circ \sigma(2,2) \circ  
(\gen21 \ot \gen21),&
\end{eqnarray*}
where $\sigma(2,2)\in \Sigma_4$ is the permutation 
\begin{equation}
\label{novy_tabak}
\sigma(2,2) = 
\left(
\begin{array}{cccc}
1 & 2 & 3 & 4
\\
1 & 3 & 2 & 4
\end{array}
\right)
\end{equation}
or diagrammatically
\[
\sigma(2,2) = \hskip 2mm
{
\unitlength=1.000000pt
\begin{picture}(30.00,30.00)(0.00,12.00)
\put(30.00,0.00){\makebox(0.00,0.00){\scriptsize$\bullet$}}
\put(30.00,30.00){\makebox(0.00,0.00){\scriptsize$\bullet$}}
\put(20.00,0.00){\makebox(0.00,0.00){\scriptsize$\bullet$}}
\put(20.00,30.00){\makebox(0.00,0.00){\scriptsize$\bullet$}}
\put(10.00,0.00){\makebox(0.00,0.00){\scriptsize$\bullet$}}
\put(10.00,30.00){\makebox(0.00,0.00){\scriptsize$\bullet$}}
\put(0.00,0.00){\makebox(0.00,0.00){\scriptsize$\bullet$}}
\put(0.00,30.00){\makebox(0.00,0.00){\scriptsize$\bullet$}}
\put(30.00,20.00){\line(0,-1){10.00}}
\put(20.00,20.00){\line(-1,-1){10.00}}
\put(10.00,20.00){\line(1,-1){10.00}}
\put(0.00,20.00){\line(0,-1){10.00}}
\end{picture}}\hskip .3em .
\]
\vglue .6em

\noindent
We will use the similar notation for elements of free \PROP{s}
throughout the paper. All our `flow diagrams' should be read from the
bottom to the top.

\begin{remark}
\label{zase_se_mi_o_ni_zdalo}
{\rm Enriquez and Etingof described in~\cite{enriquez-etingof:03} a
basis of the $\bfk$-linear space $\sfB(m,n)$ for arbitrary $m,n \geq 1$ as
follows.  Let $\gen 12 \in \sfB(1,2)$ be the equivalence class, in
$\sfB = \Gamma(\gen12,\gen21)/{\sf I}_{\sf B}$, of the generator $\gen
12 \in \Gamma(\gen 21,\gen 12)(1,2)$ (we use the same symbol both for
a generator and its equivalence class). Define ${\gen 12}^{[1]} := \id
\in \sfB(1,1)$ and, for $a \geq 2$, let
\[
{\gen 12}^{[a]} := 
   \gen 12 (\gen 12 \ot \id) (\gen 12 \ot \otexp {\id}2) \cdots  
   (\gen 12 \ot \otexp {\id}{(a-2)}) \in \sfB(1,a).
\]
Let ${\gen 21}_{[b]} \in \sfB(b,1)$ has the obvious similar
meaning. According to~\cite[Proposition~6.2]{enriquez-etingof:03}, the
elements
\[
({\gen 12}^{[a_1]} \ot \cdots \ot {\gen 12}^{[a_m]}) 
\circ \sigma \circ
({\gen 21}_{[b^1]} \ot \cdots \ot {\gen 21}_{[b^n]}),
\]
where $\sigma \in \Sigma_N$ for some $N \geq 1$, and $a_1 + \cdots +
a_m = b^1 + \cdots + b^m = N$, form a $\bfk$-linear basis
of~$\sfB(m,n)$.
This result can also be found in~\cite{kontsevich-soibelman:hopf-algebras}.
}
\end{remark}

We have already observed that \PROP{s}, and namely free ones, are
extremely huge objects. For instance, the space
$\Gamma(\gen12,\gen21)(m,n)$ is infinite-dimensional for any $m, n$,
and even its quotient $\sfB(m,n)$ is infinite-dimensional, as follows
from Proposition~6.2 of~\cite{enriquez-etingof:03} recalled in
Remark~\ref{zase_se_mi_o_ni_zdalo}. Therefore it might come as a
surprise that there are {\em three\/} natural gradings of
$\Gamma(\gen12,\gen21)(m,n)$ by finite-dimensional pieces.

Since elements of free \PROP{s} are represented by formal sums of
graphs~\cite[Section~2]{mv}, it makes sense to define the {\em
genus\/} $\genus(X)$ of a monomial $X$ in a free \PROP\ as the genus
$\dim H^1(G_X; {\bf Q})$ of the graph $G_X$ corresponding to $X$. For
example, $\genus(\ZbbZb)= \genus(\dvojiteypsilon) = 0$, while
\[
\genus(\motylek) =1.
\]

There is another grading called the {\em path grading\/} $\pth(X)$
implicitly present in~\cite{kontsevich:message}, defined as the total
number of directed paths connecting inputs with outputs of
$G_X$. Properties of the genus and path gradings are discussed
in~\cite[Section~5]{mv}. The following proposition follows immediately
from the results of~\cite{mv}.

\begin{proposition}
\label{opet_predposledni}
For any fixed $d$, the subspaces
\[
\Span\{X \in \Gamma(\gen12,\gen21)(m,n); \ \genus(X) = d\} \mbox { and }\
\Span\{X \in \Gamma(\gen12,\gen21)(m,n); \ \pth(X) = d\}
\]
are finite dimensional.
\end{proposition}

The following formula relating the path and genus gradings
was also derived in~\cite{mv}:
\begin{equation}
\label{netekla_voda}
\pth(X) \leq mn(\genus(X) +1)\ \mbox { for $X \in
\Gamma(\gen12,\gen21)(m,n)$}.
\end{equation} 

There is, of course, also the {\em obvious grading\/} $\grd(X)$ given by the
number of vertices of the graph $G_X$. Using this grading, the
decomposables of a free \PROP\ can be described as
\[
D(\Gamma(E)) = \Span\{X \in \Gamma(E);\ \grd(X) \geq 2\}.
\]
Let us recall the following important
definition~\cite{kontsevich:message,mv}.

\begin{definition}
\label{def:small}
A \hPROP\ is a collection $\sfs =\{\sfs(m,n)\}$ of dg
$(\Sigma_m,\Sigma_n)$-bimodules $\sfs(m,n)$ defined for all couples
of natural numbers except $(m,n) = (1,1)$, together with compositions
\begin{equation}
\label{circi}
\circ_i:   \sfs(m_1,n_1) \ot \sfs(1,l) 
\to \sfs(m_1,n_1+l-1),\ 1 \leq i \leq n_1,
\end{equation}
and
\begin{equation}
\label{jcirc}
\jcirc j : \sfs(k,1) \ot \sfs(m_2,n_2) 
\to \sfs(m_2 + k -1,n_2),\ 1 \leq j \leq m_2,
\end{equation}
that satisfy the axioms satisfied by operations $\circ_i$ and $\jcirc
j$, see~(\ref{Zbraslavice1}),~(\ref{Zbraslavice2}), in a general
\PROP.
\end{definition}

\begin{remark}
{\rm Observe that \hPROP{s} as introduced above cannot have a unit $1
\in \sfs(1,1)$.  We choose this convention from the following reasons.
There exist an obvious unital version of \hPROP{s}, but for all
examples of interest, including $\frac12$bialgebras, the corresponding
unital \hPROP\ would satisfy $\sfs(1,1) \cong \bfk$.
Since there clearly exists a canonical one-to-one correspondence
between unital \hPROP{s} enjoying this property and non-unital
\hPROP{s} in the sense of the above definition, the unit
would carry no information. 

Moreover, working without units enables one to define the `obvious
grading' $\grd(-)$ of free \hPROP{s} in a very natural way, without using
graphs. The same reason lead us in~\cite{markl:zebrulka} to introduce
pseudo-operads as non-unital versions of operads. The above
considerations do not apply to \PROP{s} because $\sfP(1,1)$ is
typically an infinite-dimensional space.  }
\end{remark}

Let us denote by $\Gamma_{\frac12}(\gen12,\gen21)$ the free \hPROP\
generated by operations $\gen12$ and $\gen21$.  The following
proposition, which follows again from~\cite{mv}, gives a
characterization of the subspaces
\[
\Gamma_{\frac12}(\gen12,\gen21)(m,n) \subset \Gamma(\gen12,\gen21)(m,n)
\]
in terms of the genus and path gradings introduced above.

\begin{proposition}
\label{Zbraslavice_2003}
The subspace $\Gamma_{\frac12}(\gen12,\gen21)(m,n)$ is, for $(m,n)
\not= (1,1)$, spanned by all monomials $X \in
\Gamma(\gen12,\gen21)(m,n)$ such that (i) $\genus(X) = 0$ and (ii)
$\pth(X) = mn$.  

Equivalently, $\Gamma_{\frac12}(\gen12,\gen21)(m,n)$
is the span of elements of the form $U \circ V$ with some monomials $U
\in \Gamma(\gen21)(m,1)$ and $V \in \Gamma(\gen 12)(1,n)$.
\end{proposition}

Loosely speaking, elements of $\Gamma_{\frac12}(\gen12,\gen21)$ are
formal sums of graphs made of two trees grafted by their roots.
Now it is completely obvious 
that $\Gamma_{\frac12}(\gen12,\gen21)(m,n)$ is
finite-dimensional for any $m$ and $n$. The following example shows
that both assumptions (i) and (ii) in
Proposition~\ref{Zbraslavice_2003} are necessary.

\begin{example}
\label{v_Ratajich}
{\rm
It is clear that 
$\genus(\levaplastev) = 0$,
$\pth(\levaplastev)$ = 3, and it is indeed almost obvious that
$\levaplastev \not\in \Gamma_{\frac12}(\gen12,\gen21)(2,2)$. An example
for which (ii) is satisfied but (i) is violated is provided by
\[
\genus(\motylekmeziradky) = 1 \mbox { and } \pth(\motylekmeziradky) = 4.
\]
Proposition~\ref{Zbraslavice_2003} then gives a rigorous proof of the
more or less obvious fact that
\begin{equation}
  \label{eq:4}
\motylekmeziradky \not\in \Gamma_{\frac12}(\gen12,\gen21)(2,2).
\end{equation}
On the other hand, $\genus(\dvojiteypsilon) = 0$ and 
$\pth(\dvojiteypsilon) = 4$, which corroborates that
$\dvojiteypsilon \in \Gamma_{\frac12}(\gen12,\gen21)(2,2)$.  
}
\end{example}

\begin{observation}
Bialgebras cannot be defined over \hPROP{s},
because the compatibility axiom~(\ref{compatibility}) contains an
element which does not belong to
$\Gamma_{\frac12}(\gen12,\gen21)(2,2)$, see~(\ref{eq:4}).  
\end{observation}

In contrast, $\frac12$bialgebras are algebras over the \hPROP\ $\hsfb$
defined as
\[
\hsfb := \Gamma_{\frac12}(\gen12,\gen21)/{\sf i}_{\sshsfb},
\]
with the ideal ${\sf i}_{\sshsfb}$ generated by
\[
\ZbbZb - \bZbbZ,\
\ZvvZv - \vZvvZ\ \mbox { and }\ \dvojiteypsilon.
\]

For a generator $\xi$ of biarity $(m,n)$, let $\sigmaspan(\xi) : =
\bfk[\Sigma_m] \otimes \bfk \cdot \xi \otimes \bfk[\Sigma_n]$, with
the obvious mutually compatible left $\Sigma_m$- right
$\Sigma_n$-actions.  

The first step in pursuing the Principle formulated in
Section~\ref{intro} is to describe a minimal model of the \hPROP\
$\hsfb$ for $\frac12$bialgebras in the category of \hPROP{s}. This can be
done as follows. Theorem~18 of~\cite{mv} implies that
$\hsfb$ is a Koszul quadratic \hPROP, therefore its minimal model is
given by the cobar dual $\Omega_{\frac12\tt P}(\hsfb^!)$ of the
quadratic dual $\hsfb^!$ of $\hsfb$. This cobar dual is, by definition, a
dg-\hPROP\ of the form $(\Gamma_{\frac12}(\Xi), \pa_0)$, with
\[
\Xi := \Lambda \downarrow (\hsfb^!)^*,
\]
where $\Lambda$ denotes the sheared suspension~\cite{gan},
$\downarrow$ the usual desuspension of a graded vector space and $(-)^*$
the linear dual. Because, by~\cite[Example~16]{mv}, $\hsfb^!(m,n)
\cong \bfk$ for any $(m,n) \not= (1,1)$, one immediately sees that
$\Xi := \sigmaspan(\{\xi^m_n\}_{m,n \in I})$ with
\[
I := \{ m,n \geq 1,\ (m,n) \not= (1,1)\}.
\] 
where the generator $\xi^m_n$ of biarity $(m,n)$ has degree $n+m-3$.

It remains to describe the differential $\pa_0$ which is, 
by definition, the unique derivation extending the linear dual of the
structure operations of $\hsfb^!$. The result is given in the
following theorem.
 
\begin{theorem}
\label{jsem_prepracovany}
There is a minimal model of the \hPROP\ $\hsfb$
\begin{equation}
\label{eq1}
(\hsfb, \partial = 0) \stackrel{\rho_{\frac 12}}{\longleftarrow}
(\Gamma_{\frac12}(\Xi), \pa_0),
\end{equation}
with the map $\rho_{\frac 12}$ defined by
\[
\rho_{\frac 12}(\xi^1_2) := \jednadva, \
\rho_{\frac 12}(\xi^2_1) := \dvajedna,
\]
while $\rho_{\frac 12}$ is trivial on all 
remaining generators. The differential
$\pa_0$ is given by the formula
\begin{equation}
\label{Paal_v_Praze}
\pa_0 (\xi^m_n) :=
(-1)^{m} \xi_1^m \circ \xi^1_n +
\sum_U (-1)^{i(s+1) + m} \xi_u^m \circ_i \xi^1_s 
+\sum_V (-1)^{j(t+1) + 1} \xi_1^t \hskip 2pt \cric j \xi^v_n,
\end{equation}
where we set $\xi^1_1 := 0$,
\[
U := \{u,s \geq 1, \ u+s =n +1, \ 1 \leq i \leq u
\}
\]
and
\[
V = \{t,v \geq 1, \ t+v = m + 1, \ 1 \leq j \leq v
\}.
\]
\end{theorem}

It follows from the remarks preceding Theorem~\ref{jsem_prepracovany}
that a quadratic Koszul \hPROP\ admits a canonical functorial minimal
model, given by the cobar dual of its quadratic dual. It can also be
proved that minimal models of \hPROP{s} are unique up to
isomorphism.

\begin{example}
\label{Orlik}
{\rm
If we denote $\xi^1_2 = \jednadva$ and $\xi^2_1 = \dvajedna$, then
$\pa_0(\jednadva) = \pa_0(\dvajedna) =0$.
If $\xi^2_2 = \dvadva$, then
\[
\pa_0(\dvadva) = \dvojiteypsilon. 
\]
With the obvious, similar notation,
\begin{eqnarray}
\label{jedna-tri}
\pa_0(\jednatri) &=& \ZbbZb - \bZbbZ,
\\
\label{jedna-ctyri}
\pa_0(\jednactyri) &=& \ZbbZbb - \bZbbZb + \bbZbbZ - \ZbbbZb - \bZbbbZ,
\\
\nonumber 
\pa_0(\trijedna) &=& \ZvvZv - \vZvvZ,
\\
\nonumber 
\pa_0(\dvatri) &=& 
\dvacarkatri - \dvaZbbZb + \dvabZbbZ,
\\
\nonumber 
\pa_0(\tridva) &=&  -
{
\unitlength=.2pt
\begin{picture}(40.00,50.00)(0.00,-50.00)
\put(20.00,-30.00){\line(0,1){10.00}}
\put(20.00,0.00){\line(0,-1){20}}
\bezier{20}(20.00,-20.00)(30.00,-10.00)(40.00,0.00)
\bezier{20}(20.00,-20.00)(10.00,-10.00)(0.00,0.00)
\bezier{20}(20.00,-30.00)(30.00,-40.00)(40.00,-50.00)
\bezier{20}(0.00,-50.00)(10.00,-40.00)(20.00,-30.00)
\end{picture}}
+ \ZvvZvdva - \vZvvZdva,
\\
\nonumber 
\pa_0(\tritri)&=& 
-
{
\unitlength=.4pt
\begin{picture}(24.00,30.00)(-2.00,0.00)
\bezier{20}(10.00,10.00)(15.00,5.00)(20.00,0.00)
\bezier{20}(0.00,0.00)(5.00,5.00)(10.00,10.00)
\put(0,-5){
\bezier{20}(10.00,20.00)(15.00,25.00)(20.00,30.00)
\bezier{20}(0.00,30.00)(5.00,25.00)(10.00,20.00)
\put(10.00,30.00){\line(0,-1){25}}
}
\end{picture}}
+ \triZbbZb - \tribZbbZ + \ZvvZvtri - \vZvvZtri,
\\
\nonumber 
\pa_0(\dvactyri) &=& 
{
\unitlength=.1pt
\begin{picture}(96.00,120.00)(-8.00,0.00)
\qbezier(40.00,60.00)(20.00,80.00)(0.00,100.00)
\qbezier(40.00,60.00)(60.00,80.00)(80.00,100.00)
\put(40.00,60.00){\line(0,-1){20.00}}
\qbezier(40.50,39.00)(49.00,14.50)(53.50,0.00)
\qbezier(40.00,40.00)(32.50,16.50)(26.00,0.00)
\qbezier(40.00,40.00)(60.00,20.00)(80.00,0.00)
\qbezier(40.00,40.00)(20.00,20.00)(0.00,0.00)
\end{picture}}
+ \dvaZbbbZb + \dvabZbbbZ - \dvaZbbZbb + \dvabZbbZb - \dvabbZbbZ,~\mbox{etc.}
\end{eqnarray}
Observe that~(\ref{Paal_v_Praze}) for $m = 1$ gives
\[
\pa_0(\xi^1_n) = 
\sum_U (-1)^{i(s+1) +1} \xi_u^1 \circ_i \xi^1_s,
\]
where $U$ is as in Theorem~\ref{jsem_prepracovany}.
Therefore the sub-\hPROP\ generated by $\xi^1_2$, $\xi^1_3$,
$\xi^1_4, \ldots$ is in fact isomorphic to
the minimal model ${\mathcal A}_\infty$ for the operad of
associative algebras as described in~\cite{markl:zebrulka}.

It is well-known that ${\mathcal A}_\infty$ is the operad of
cellular chains of a cellular topological operad $\underline {\mathcal K}
= \{K_n\}_{n \geq 2}$ such that each $K_n$ is an $(n-2)$-dimensional
convex polyhedron -- the Stasheff associahedron 
(see~\cite[Section~1.6]{markl-shnider-stasheff:book}). The formulas for
the differential $\pa_0(\xi^1_n)$ then reflect the 
decomposition of the topological boundary of the top dimensional cell
of $K_n$ into the union of codimension one faces. For example, the two
terms in the right-hand side of~(\ref{jedna-tri}) correspond to the two
endpoints of the interval $K_3$, the five terms in the right-hand side
of~(\ref{jedna-ctyri}) to the five edges of the pentagon $K_4$,~etc.
}
\end{example}

\begin{remark}
\label{kasparek}
{\rm 
Just as there are non-$\Sigma$ operads as simplified versions of
operads without the actions of symmetric
groups~\cite[Definition~II.1.14]{markl-shnider-stasheff:book}, there
are obvious notions of non-$\Sigma$ \PROP{s} and non-$\Sigma$
\hPROP{s}. Of a particular importance for us will be the free
non-$\Sigma$ \hPROP\ $\uGamma_{\frac12}(\uXi)$ generated by $\uXi :=
\span(\{\xi^m_n\}_{m,n \in I})$, where $\xi^m_n$ and $I$ are as in
Theorem~\ref{jsem_prepracovany}. There clearly exists, for any $m$ and
$n$, a $\pa_0$-invariant factorization of $\Sigma_m$-$\Sigma_n$ spaces
\begin{equation}
\label{pristi_tyden_jedu_za_Jitkou}
\Gamma_{\frac 12}(\Xi)(m,n) \cong \bfk[\Sigma_m] \otimes
\uGamma_{\frac12}(\uXi)(m,n) \otimes \bfk[\Sigma_n].
\end{equation}

Therefore, the acyclicity of $(\Gamma_{\frac 12}(\Xi),\pa_0)$ is
equivalent to the acyclicity of
$(\uGamma_{\frac12}(\uXi),\pa_0)$. Observe that there is no analog of
factorization~(\ref{pristi_tyden_jedu_za_Jitkou}) for \PROP{s}.
}
\end{remark}

\begin{remark}
{\em Another way to control the combinatorial explosion of \PROP{s}
was suggested by W.L.~Gan who introduced {\em dioperads\/}.
Roughly speaking, a dioperad is a \PROP\ in which only compositions
based on graphs of genus zero are allowed, see~\cite{gan} for
details.

Dioperads are slightly bigger than \hPROP{s}. The piece
$\Gamma_D(\gen12,\gen21)(m,n)$ of the free dioperad
$\Gamma_D(\gen12,\gen21)$ is spanned by genus zero monomials of
$\Gamma(\gen12,\gen21)(m,n)$, with no restriction on the path
grading. Therefore, for instance,
\[
\levaplastev \in \Gamma_D(\gen12,\gen21)(2,2), \mbox { while } 
\levaplastev \not\in \Gamma_{\frac12}(\gen12,\gen21)(2,2),
\]
see Example~\ref{v_Ratajich}. The relation between \PROP{s},
dioperads and \hPROP{s} is analyzed in~\cite{mv}, where we also
explain why \hPROP{s} are better suited for our purposes than dioperads. 
 
}
\end{remark}

Let us finish Step~1 formulated in Section~\ref{intro} by describing a
minimal model of the \PROP\ $\sfhB$, following again~\cite{mv}.
Observe first that the \PROP\ $\sfhB$ is generated by the \hPROP\
$\sfhb$. By this we mean that $\sfhB = L_{\ssBox}(\sfhb)$, where
$L_{\ssBox} : \cathPROP \to \catPROP$ is the left adjoint to the
forgetful functor $\Box :\catPROP \to \cathPROP$.  The functor
$L_{\ssBox}$ is, by~\cite[Theorem~4]{mv}, {\em exact\/}.  This
surprisingly deep statement follows from the fact, observed by
M.~Kontsevich in~\cite{kontsevich:message}, that $L_{\ssBox}$ is a
{\em polynomial\/} functor in the sense recalled
in~\cite[Definition~1]{kontsevich-soibelman:hopf-algebras}. The last
thing we need to realize is that $L_{\ssBox}(\Gamma_{\frac12}(\Xi),
\pa_0) = (\Gamma(\Xi), \pa_0)$, where the differential $\pa_0$ is in
both cases given by the same formula on the space of generators. We
conclude that the application of the functor $L_{\ssBox}$ to the
minimal model of the \hPROP\ $\sfhb$ described in
Theorem~\ref{jsem_prepracovany} gives a minimal model of the \PROP\
$\sfhB$. We obtain

\begin{theorem}
\label{dnes_tatarak_u_Jany}
The dg-\PROP\
\begin{equation}
{\sf M}_0 := (\Gamma(\Xi),\pa_0),
\end{equation}
where the generators $\Xi$ are as in Theorem~\ref{jsem_prepracovany}
and the differential $\pa_0$ is given by formula~(\ref{Paal_v_Praze}),
is a minimal model of the \PROP\ $\sfhB$ for
$\frac12$bialgebras.
\end{theorem}

\begin{remark}
{\rm\
For a \hPROP\ $\sfs$, let $P(\sfs)$ be the augmented \PROP\ 
whose augmentation ideal equals $\sfs$,
whose compositions $\circ_i$ and $\cric j$ of~(\ref{Zbraslavice1}) 
and~(\ref{Zbraslavice2}) are those of $\sfs$, and other compositions
(that is, those not allowed for \hPROP{s}) are set to be zero. 
Theorem~\ref{dnes_tatarak_u_Jany} expresses the fact
that the \PROP\ $P(\hsfb^!)$ is the quadratic dual of the \PROP\ $\sfhB$
in the category of \PROP{s} in the sense of 
B.~Vallette~\cite{vallette,vallette:CMR04}. 
}
\end{remark}

\section{Main theorem and the proof - first attempt}
\label{zase_starosti_mam}

Let us formulate the main theorem of the paper.

\begin{theorem}
\label{main}
There exists a minimal model $({\sf M},\pa)$ of the \PROP\ \sfB\ for
bialgebras that is a perturbation of the minimal model $({\sf
M}_0,\pa_0)$ of the \PROP\ $\sfhB$ for $\frac12$bialgebras described
in Theorem~\ref{dnes_tatarak_u_Jany}. By this
we mean that
\[
({\sf M},\pa)  = (\Gamma(\Xi),\pa_0 + \papert), 
\]
where the generators $\Xi =  \sigmaspan(\{\xi^n_m\}_{m,n \in I})$ are
as in Theorem~\ref{jsem_prepracovany} and $\pa_0$ is a derivation given by
formula~(\ref{Paal_v_Praze}). The perturbation $\papert$ raises the
genus and preserves the path grading. More precisely, 
$\papert = \pa_1 + \pa_2 + \pa_3 + \cdots$,
where $\pa_g$ raises the genus by $g$, preserves the path grading and,
moreover, 
\begin{equation}
\label{L}
\pa_g(\xi^m_n) = 0 \mbox { for } g > (m-1)(n-1).
\end{equation}
\end{theorem}

Uniqueness of minimal models for \PROP{s} is discussed in
Section~\ref{po_zavodech_v_Hosine}.  Observe that~(\ref{L}) implies
$\pa(\xi^1_n) = \pa_0(\xi^1_n)$ for all $n$. Therefore the
sub-dg-operad generated in $({\sf M},\pa)$ by
$\xi^1_2,\xi^1_3,\xi_4^1,\ldots$ is isomorphic to the operad
describing strongly homotopy associative algebras.

Formulas for the perturbed differential $\papert(\xi^m_n)$ are, for
some small $m$ and $n$, given in Section~\ref{explicit}. Although
Theorem~\ref{main} does not describe the perturbation $\papert$
explicitly, it describes the space of generators $\Xi$ of the
underlying free \PROP.  This itself seems to be very nontrivial
information. It will also be clear later that $\pa_0$ is in fact the
quadratic part (with respect to the `obvious' grading recalled in
Section~\ref{zase_starosti_mam}) of the perturbed differential $\pa$,
therefore, using the terminology borrowed from rational homotopy
theory, the unperturbed model $(\sfM_0,\pa_0)$ describes the `homotopy
Lie algebra' of the \PROP\ $\sfB$.

Let us try to prove Theorem~\ref{main} by
constructing na\"{\i}vely a perturbation $\papert$ as
\[
\papert = \pa_1 + \pa_2 + \pa_3 + \cdots,
\]
where each $\pa_g$ is a derivation raising the genus by $g$. Observe that
$\pa_g(\xi^m_n)$ must be a sum of decomposable elements, because the
generators are of genus $0$. It is, of course, enough to define
$\papert$ on the generators $\xi^m_n \in \Xi$ and extend it as a
derivation.

We construct $\papert(\xi^m_n)$ inductively. Let $N:= m+n$.
For $N =3$, we must put
\[
\papert(\gen12) = \papert(\gen21) = 0. 
\]
Also for $N = 4$ the formula for the differential is dictated by the
axioms of bialgebras:
\[
\papert(\gen31) := \papert(\gen13) := 0 
\mbox { and } \papert(\gen22) := -\motylek.
\] 
For $N=5$ we put 
\[
\papert(\gen 14) = \papert(\gen 41) :=0;
\]
$\papert(\gen 23)$ and $\papert(\gen 32)$ are given by
formulas
\begin{eqnarray}
\label{A} \hskip 2mm
\papert(\gen 23)
\hskip -2mm &:=& \hskip -2mm (\gen12 \otm \gen 12) \hskip -.26em \circ
\hskip -.26em\sigma(2,2) \hskip -.26em \circ \hskip -.26em (\gen 21 \otm
\gen 22 - \gen 22 \otm \gen 21) - (\bZbbZ \otm \gen 13 + \gen 13 \otm
\ZbbZb) \hskip -.26em \circ \hskip -.26em\sigma(3,2)\hskip -.26em \circ
\hskip -.26em (\gen 21 \otm \gen 21 \otm\gen 21), 
\\
\label{BB} \hskip 2mm
\papert(\gen 32) 
\hskip -2mm &:=& \hskip -2mm (\gen 12 \otm \gen 22 - \gen 12 \otm \gen
22) \hskip -.26em \circ \hskip -.26em\sigma(2,2)\hskip -.26em \circ
\hskip -.26em (\gen 21 \otm \gen 21) + (\gen 12 \otm \gen 12 \otm \gen
12) \hskip -.26em \circ \hskip -.26em\sigma(2,3)\hskip -.26em \circ
\hskip -.26em (\vZvvZ \otm \gen 31 + \gen 31 \otm \ZvvZv).
\end{eqnarray}
In the above displays, $\sigma(2,2)$ is the same as
in~(\ref{novy_tabak}), 
\begin{equation}
\label{prijede_Jim}
\sigma(3,2) := 
\left(
\begin{array}{cccccc}
1 & 2 & 3 & 4 & 5 & 6
\\
1 & 4 & 2 & 5 & 3 & 6
\end{array}
\right)
= \hskip 3mm
{
\unitlength=1.000000pt
\begin{picture}(50.00,30.00)(0.00,12.00)
\put(50.00,0.00){\makebox(0.00,0.00){\scriptsize$\bullet$}}
\put(50.00,30.00){\makebox(0.00,0.00){\scriptsize$\bullet$}}
\put(40.00,0.00){\makebox(0.00,0.00){\scriptsize$\bullet$}}
\put(40.00,30.00){\makebox(0.00,0.00){\scriptsize$\bullet$}}
\put(30.00,0.00){\makebox(0.00,0.00){\scriptsize$\bullet$}}
\put(30.00,30.00){\makebox(0.00,0.00){\scriptsize$\bullet$}}
\put(20.00,0.00){\makebox(0.00,0.00){\scriptsize$\bullet$}}
\put(20.00,30.00){\makebox(0.00,0.00){\scriptsize$\bullet$}}
\put(10.00,0.00){\makebox(0.00,0.00){\scriptsize$\bullet$}}
\put(10.00,30.00){\makebox(0.00,0.00){\scriptsize$\bullet$}}
\put(0.00,0.00){\makebox(0.00,0.00){\scriptsize$\bullet$}}
\put(0.00,30.00){\makebox(0.00,0.00){\scriptsize$\bullet$}}
\put(0.00,30.00){\line(0,1){0.00}}
\put(50.00,20.00){\line(0,-1){10.00}}
\put(40.00,20.00){\line(-1,-1){10.00}}
\put(30.00,20.00){\line(-2,-1){20.00}}
\put(20.00,20.00){\line(2,-1){20.00}}
\put(10.00,20.00){\line(1,-1){10.00}}
\put(0.00,20.00){\line(0,-1){10.00}}
\end{picture}} \hskip 2mm
\end{equation}

\noindent 
with our usual convention that the `flow diagrams'
should be read from the bottom to the top, and $\sigma(3,2) :=
\sigma(2,3)^{-1}$. Higher terms of the perturbed differential can be
constructed by the standard homological perturbation theory as
follows.

Suppose we have already constructed $\papert(\xi^u_v)$ for all $u+v <
N$ and fix some $m$ and $n$ such that 
$m+n = N > 5$. We are looking for $\papert(\xi^m_n)$
of the form
\begin{equation}
\label{uplne_hotovy}
\papert(\xi^m_n) = \pa_1(\xi^m_n) + \pa_2(\xi^m_n) + \pa_3(\xi^m_n) + \cdots
\end{equation}
where $\genus(\pa_g(\xi^m_n)) = g$.
Condition $(\pa_0 + \papert)^2(\xi^m_n) = 0$ can be rewritten as
\[
\sum_{s+t = g} \pa_s\pa_t(\xi^m_n) = 0 \mbox { for each $g \geq 1$.}
\]
We must therefore find inductively elements $\pa_g(\xi^m_n)$, $g \geq
1$, solving the equation
\begin{equation}
\label{budu_na_tech_zavodech_litat?}
\pa_0\pa_g(\xi^m_n) = 
-\sum_\doubless{s+t = g}{t < g} \pa_s\pa_t(\xi^m_n).
\end{equation}

Observe that the right-hand side
of~(\ref{budu_na_tech_zavodech_litat?}) makes sense, because
$\pa_t(\xi^m_n)$ is a combination of $\xi^u_v$'s with $u+v <
N$, therefore $\pa_s\pa_t(\xi^m_n)$ has already been defined. To
verify that the right-hand side of~(\ref{budu_na_tech_zavodech_litat?})
is a $\pa_0$-cycle is also easy:
\begin{eqnarray*}
\pa_0(-\sum_\doubless{s+t = g}{t < g}  \pa_s\pa_t(\xi^m_n))
&=& -\sum_\doubless{s+t = g}{t < g} \pa_0\pa_s\pa_t(\xi^m_n)
=
\sum_\doubless{s+t = g}{t < g}\sum_\doubless{a+b=s}{b < s}
\pa_a\pa_b\pa_t(\xi^m_n)
\\
&=&
\sum_{1 \leq i \leq g}\pa_i (\sum_{k+l = g-i} \pa_k\pa_l(\xi^m_n)) =0.
\end{eqnarray*}

The degree of the right-hand side
of~(\ref{budu_na_tech_zavodech_litat?}) is $N-5$, which is a positive
number, by our assumption $N > 5$. This implies 
that~(\ref{budu_na_tech_zavodech_litat?}) has a solution, 
because $(\Gamma(\Xi), \pa_0)$ is, by Theorem~\ref{dnes_tatarak_u_Jany}, 
$\pa_0$-acyclic in positive dimensions.%
\qed

There is however a serious flaw in the above proof: there is \uu{no}
\uu{reason} \uu{to} \uu{assume} \uu{that} \uu{the}
\uu{sum}~(\ref{uplne_hotovy}) \uu{is} \uu{finite}, that is, that the
right-hand side of~(\ref{budu_na_tech_zavodech_litat?}) is trivial for
$g$ sufficiently large!!! This convergence problem can be fixed by
finding subspaces $F(m,n) \subset \Gamma(\Xi)(m,n)$ satisfying the
properties listed in the following definition.

\begin{definition}
\label{s_tou_Elisabetou_jsem_to_prepisk}
The collection $F$ of subspaces  $F(m,n) \subset
\Gamma(\Xi)(m,n)$ is {\em friendly\/} if
\begin{itemize}
\item[(i)]
for each $m$ and $n$, there exists a constant $C_{m,n}$ such that $F(m,n)$ 
does not contain elements of genus $> C_{m,n}$,
\item[(ii)] $F$ is stable under all derivations (not
necessary differentials) $\omega$ satisfying $\omega(\Xi) \subset F$,
\item[(iii)]
\vskip -2,5mm
$\pa_0(\Xi) \subset F$, $\motylek \in F(2,2)$, the right-hand side 
of~(\ref{A}) belongs to $F(2,3)$ and the right-hand side 
of~(\ref{BB}) belongs to $F(3,2)$, and
\item[(iv)]
$F$ is $\pa_0$-acyclic in positive degrees. 
\end{itemize}
\end{definition}

Observe that (ii) with (iii) imply that $F$ is $\pa_0$-stable,
therefore (iv) makes sense. Observe also that we do not demand
$F(m,n)$ to be $\Sigma_m$-$\Sigma_n$ invariant.

Suppose we are given such a friendly collection.
We may then, in the above na{\"\i}ve proof, assume inductively that
\begin{equation}
\label{jim_spendliky}
\papert(\xi^m_n) \in F(m,n).
\end{equation}
Indeed,~(\ref{jim_spendliky}) is satisfied for $m+n = 3,4,5$, by (iii).
Condition~(ii) guarantees that the right-hand side
of~(\ref{budu_na_tech_zavodech_litat?}) belongs to $F(m,n)$,
while~(iv) implies that~(\ref{budu_na_tech_zavodech_litat?}) can be
solved in $F(m,n)$. Finally, (i) guarantees, in the
obvious way, the convergence.
 
In this paper, we use the friendly collection $\sfS \subset
\Gamma(\Xi)$ of {\em special elements\/}, introduced in
Section~\ref{special}. The collection $\sfS$ is generated by the free
non-$\Sigma$ \hPROP\ $\uGamma_{\frac12}(\uXi)$, see
Remark~\ref{kasparek}, by a suitably
restricted class of compositions that naturally generalize those
involved in $\motylekmeziradky$.

Another possible choice was proposed in~\cite{mv}, namely the friendly
collection defined by
\[
F(m,n) := \{f \in \Gamma(\Xi);\ \pth(f) = mn\}.
\]
This choice is substantially bigger than the collection of special
elements and contains `strange' elements, such as
\begin{center}
{
\unitlength=.2pt
\begin{picture}(110.00,120.00)(0.00,0.00)
\put(110.00,60.00){\makebox(0.00,0.00)[l]{$\in F(2,2)$}}
\put(20.00,60.00){\line(0,-1){60.00}}
\qbezier(80.00,40.00)(70.00,30.00)(60.00,20.00)
\qbezier(60.00,60.00)(70.00,50.00)(80.00,40.00)
\qbezier(40.00,40.00)(50.00,30.00)(60.00,20.00)
\qbezier(60.00,60.00)(50.00,50.00)(40.00,40.00)
\qbezier(40.00,100.00)(50.00,90.00)(60.00,80.00)
\qbezier(40.00,100.00)(30.00,90.00)(20.00,80.00)
\qbezier(0.00,100.00)(10.00,90.00)(20.00,80.00)
\put(60.00,20.00){\line(0,-1){20.00}}
\put(60.00,80.00){\line(0,-1){20.00}}
\put(20.00,80.00){\line(0,-1){20.00}}
\put(40.00,100.00){\line(0,1){20.00}}
\put(0.00,120.00){\line(0,-1){20.00}}
\end{picture}}
\end{center}
which we certainly do not want to consider. We believe that special
elements are, in a suitable sense, the smallest possible friendly collection.
 
Properties of special elements are studied in Sections~\ref{kam_nas_vystehuji?}
and~\ref{ac}. Section~\ref{proof_problems} then contains a proof of
Theorem~\ref{main}.

\section{Special elements}
\label{special}

We introduce, in Definition~\ref{6}, special elements in arbitrary
free \PROP{s}. We need first the following:

\begin{definition}
\label{1}
For $k ,l \geq 1$ and $1 \leq i \leq kl$, let $\sigma(k,l) \in
\Sigma_{kl}$ be the permutation given~by
\[
\sigma(i) := k(i-1 - (s-1)l) + s,
\]
where $s$ is such that $(s-1)l < i \leq sl$. We call permutations of
this form {\em special permutations\/}.
\end{definition}

To elucidate the nature of these permutations, suppose we
have associative algebras $U_1,\ldots,U_k$. The above permutation is
exactly the permutation used to define the induced associative
algebra structure on the product
\[
(U_1 \ot \cdots \ot U_k) \otimes \cdots  \ot (U_1 \ot \cdots \ot U_k) 
\mbox {\hskip 4mm /$l$-times/},
\]
that is, the permutation which takes
\[
(u^1_1 \ot u^1_2 \ot \cdots \ot u^1_k) \ot 
(u^2_1 \ot  u^2_2 \ot \cdots  \ot  u^2_k) 
\ot \cdots \ot (u^l_1 \ot  u^l_2 \ot\cdots  \ot u^l_k)  
\]
to
\[
(u_1^1 \ot u^2_1 \ot \cdots \ot  u^l_1) \ot
(u_2^1 \ot u^2_2 \ot \cdots  \ot u^l_2) \ot \cdots
\ot (u_k^1 \ot u^2_k \ot \cdots \ot  u^l_k). 
\]

\begin{example}
\label{2}
{\rm\ We have already seen examples of special permutations: the
permutation $\sigma(2,2)$ in~(\ref{novy_tabak}) and the permutation
$\sigma(3,2)$ in~(\ref{prijede_Jim}).  Observe that, for arbitrary
$k,l \geq 1$, $\sigma(k,1) = 1_{\Sigma_k}$, $\sigma(1,l) =
1_{\Sigma_l}$ and $\sigma(k,l) = \sigma(l,k)^{-1}$.  }
\end{example}

Special elements are defined using a special class of compositions
defined as follows.

\begin{definition}
\label{3}
Let $\sfP$ be an arbitrary \PROP. 
Let $k,l \geq 1$, $\Rada a1l \geq 1$, $\Rada b1k \geq 1$,   
$\Rada A1l \in \sfP(a_i,k)$ and $\Rada B1k \in \sfP(l,b_j)$.
Then define the {\em $(k,l)$-fraction\/}
\[
\frac{A_1 \cdots A_l}{B_1 \cdots B_k} :=
(A_1 \otimes \cdots \otimes A_l) 
\circ \sigma(k,l) \circ 
(B_1 \otimes \cdots \otimes B_k)
\in \sfP(a_1+\cdots + a_l,b_1+\cdots + b_k).
\]
\end{definition}

\begin{example}
{\rm\
If $k=1$ or $l=1$, the $(k,l)$-fractions give the `operadic'
compositions:
\[
\frac{A_1 \ot \cdots \ot A_l}{B_1} = (A_1 \ot \cdots \ot A_l) \circ B_1
\ \mbox { and }\
\frac{A_1}{B_1 \ot \cdots \ot B_k} = A_1 \circ (B_1 \ot \cdots \ot B_k).
\]
}
\end{example}

We are going to use `dummy variables,' that is, 
for instance, $A \in \sfP(*,n)$ for a fixed $n
\geq 1$ means that $A \in \sfP(m,n)$ for {\em some\/} $m \geq 1$.

\begin{example}
\label{4}
{\rm
For $\zeroatwo a, \zeroatwo b \in \sfP(*,2)$ and $\twoazero c,
\twoazero d \in \sfP(2,*)$,

\vglue -1.3em
\[
\frac{\zeroatwo a \hskip .2em  \zeroatwo b}%
     {\twoazero c \hskip .2em \twoazero d}
=
(\raisebox{-2pt}{\zeroatwo a} \hskip -.3em 
\ot \raisebox{-2pt}{\zeroatwo b} 
\hskip -4pt) 
\circ \sigma(2,2) \circ 
(\raisebox{-2pt}{\twoazero c} 
\hskip -.3em\ot \raisebox{-2pt}{\twoazero d} \hskip -4pt)
=
{
\unitlength=.5pt
\begin{picture}(70.00,70.00)(0.00,30.00)
\put(50.00,10.00){\makebox(0.00,0.00)[l]{$d$}}
\put(10.00,10.00){\makebox(0.00,0.00)[l]{$c$}}
\put(50.00,60.00){\makebox(0.00,0.00)[l]{$b$}}
\put(10.00,60.00){\makebox(0.00,0.00)[l]{$a$}}
\put(60.00,50.00){\line(0,-1){30.00}}
\put(40.00,0.00){\line(0,1){20.00}}
\put(70.00,0.00){\line(-1,0){30.00}}
\put(70.00,20.00){\line(0,-1){20.00}}
\put(40.00,20.00){\line(1,0){30.00}}
\put(40.00,70.00){\line(0,-1){20.00}}
\put(70.00,70.00){\line(-1,0){30.00}}
\put(70.00,50.00){\line(0,1){20.00}}
\put(40.00,50.00){\line(1,0){30.00}}
\put(50.00,50.00){\line(-1,-1){30.00}}
\put(20.00,50.00){\line(1,-1){30.00}}
\put(0.00,0.00){\line(0,1){20.00}}
\put(30.00,0.00){\line(-1,0){30.00}}
\put(30.00,20.00){\line(0,-1){20.00}}
\put(0.00,20.00){\line(1,0){30.00}}
\put(10.00,40.00){\line(0,-1){20.00}}
\put(10.00,50.00){\line(0,-1){10.00}}
\put(30.00,70.00){\line(-1,0){30.00}}
\put(30.00,50.00){\line(0,1){20.00}}
\put(0.00,50.00){\line(1,0){30.00}}
\put(0.00,70.00){\line(0,-1){20.00}} 
\end{picture}}\hskip 2mm.
\]
Similarly, for $\zeroathree x, \zeroathree y \in \sfP(*,3)$ and
$\twoazero z, \twoazero u, \twoazero v \in \sfP(2,*)$,

\vglue -1.3em
\[
\frac{\zeroathree x \hskip 4mm \zeroathree y}%
     {\twoazero z \hskip 1mm \twoazero u \hskip 1mm \twoazero v} 
=
(\hskip -2pt \raisebox{-2pt}{\zeroathree x}\hskip -.2em \ot\hskip -.2em
\raisebox{-2pt}{\zeroathree y}\hskip -2pt) 
\circ \sigma(3,2) \circ
(\raisebox{-2pt}{\twoazero z}\hskip -.2em 
\ot  \raisebox{-2pt}{\twoazero u}  \hskip -.2em
\ot \raisebox{-2pt}{\twoazero v}\hskip -2pt) 
=
{
\unitlength=.5pt
\begin{picture}(110.00,70.00)(0.00,30.00)
\put(90.00,10.00){\makebox(0.00,0.00)[l]{$v$}}
\put(50.00,10.00){\makebox(0.00,0.00)[l]{$u$}}
\put(10.00,10.00){\makebox(0.00,0.00)[l]{$z$}}
\put(90.00,60.00){\makebox(0.00,0.00){$y$}}
\put(20.00,60.00){\makebox(0.00,0.00){$x$}}
\put(110.00,20.00){\line(-1,0){30.00}}
\put(110.00,0.00){\line(0,1){20.00}}
\put(80.00,0.00){\line(1,0){30.00}}
\put(80.00,20.00){\line(0,-1){20.00}}
\put(40.00,0.00){\line(0,1){20.00}}
\put(70.00,0.00){\line(-1,0){30.00}}
\put(70.00,20.00){\line(0,-1){20.00}}
\put(40.00,20.00){\line(1,0){30.00}}
\put(0.00,0.00){\line(0,1){20.00}}
\put(30.00,0.00){\line(-1,0){30.00}}
\put(30.00,20.00){\line(0,-1){20.00}}
\put(0.00,20.00){\line(1,0){30.00}}
\put(60.00,20.00){\line(1,1){30.00}}
\put(20.00,20.00){\line(2,1){60.00}}
\put(30.00,50.00){\line(2,-1){60.00}}
\put(20.00,50.00){\line(1,-1){30.00}}
\put(100.00,50.00){\line(0,-1){30.00}}
\put(70.00,70.00){\line(0,-1){20.00}}
\put(110.00,70.00){\line(-1,0){40.00}}
\put(110.00,50.00){\line(0,1){20.00}}
\put(70.00,50.00){\line(1,0){40.00}}
\put(0.00,70.00){\line(0,-1){20.00}}
\put(40.00,70.00){\line(-1,0){40.00}}
\put(40.00,50.00){\line(0,1){20.00}}
\put(0.00,50.00){\line(1,0){40.00}}
\put(10.00,50.00){\line(0,-1){30.00}}
\end{picture}}\hskip 2pt.
\]
}
\end{example}
If $\sfP$ is a dg-\PROP\, with differential $\pa$, then it easily
follows from Definition~\ref{3} that
\begin{eqnarray*}
\pa\left(
\frac{A_1 \cdots A_l}{B_1 \cdots B_k} 
\right)
&=&
\sum_{1 \leq i \leq l}
(-1)^{\deg(A_1) + \cdots +\deg(A_{i-1})}
\frac{A_1 \cdots \pa A_i \cdots A_l}{B_1 
\hskip 1em \cdots \cdots\hskip 1em 
B_k}
\\
&+&
\sum_{1 \leq j \leq k}
(-1)^{\deg(A_1) + \cdots + \deg(A_l) + \deg(B_1) + \cdots+ \deg(B_{j-1})}
\frac{A_1 \hskip 1em \cdots \cdots\hskip 1em  A_l}
{B_1 \cdots \pa B_j  \cdots B_k}.
\end{eqnarray*}

Suppose that the \PROP\ $\sfP$ is free, therefore the genus of monomials
of $\sfP$ is defined. It is clear that, under the notation of
Example~\ref{4},  
\[
\genus\left(
     \frac{\zeroatwo a \hskip .2em  \zeroatwo b}
          {\twoazero c \hskip .2em \twoazero d}
\right)
= 1 
+ \genus(\raisebox{-2pt}{\zeroatwo a}\hskip -.3em) 
+ \genus(\raisebox{-2pt}{\zeroatwo b}\hskip -.3em) 
+ \genus(\raisebox{-2pt}{\twoazero c}\hskip -.3em)
+ \genus(\raisebox{-2pt}{\twoazero d}\hskip -.3em) 
\]
and also that
\[
\genus\left(
\frac{\zeroathree x \hskip 4mm \zeroathree y}%
     {\twoazero z \hskip 1mm \twoazero u \hskip 1mm \twoazero v} 
\right) = 2
+\genus(\hskip -.3em\raisebox{-2pt}{\zeroathree x}\hskip -.3em)
+\genus(\hskip -.3em\raisebox{-2pt}{\zeroathree y}\hskip -.3em)
+\genus(\raisebox{-2pt}{\twoazero z}\hskip -.3em)
+\genus(\raisebox{-2pt}{\twoazero u}\hskip -.3em)
+\genus(\raisebox{-2pt}{\twoazero v}\hskip -.3em).
\]
The following lemma generalizes the above formulas.

\begin{lemma}
\label{5}
Let $\sfP$ be a free \PROP. Then the genus of the $(k,l)$-fraction is
given by
\[
\genus\left(\frac{A_1 \cdots A_l}{B_1 \cdots B_k}\right) =
(k-1)(l-1) + \sum_1^l \genus(A_i) + \sum_1^k \genus(B_j). 
\]
\end{lemma}

\begin{proof}
A straightforward and easy verification.%
\end{proof}

\begin{definition}
\label{6}
Let us define the collection $\sfS \subset \Gamma(\Xi)$ of {\em
special elements\/} to be the smallest collection of linear subspaces 
$\sfS(m,n) \subset \GammaXi(m,n)$ such that:
\begin{itemize}
\item[(i)]
$\id \in \sfS(1,1)$ and all generators $\xi^m_n \in \Xi$ belongs to
$\sfS$, and
\item[(ii)]
if $k,l \geq 1$ and $\Rada A1l, \Rada B1k \in \sfS$, then
\[
\frac{A_1 \cdots A_l}{B_1 \cdots B_k} \in \sfS.
\]

\end{itemize}
\end{definition}

\begin{remark}
\label{pocasi_nevychazi}
{\rm
One may introduce
{\em special\/} \PROP{s} as objects similar to \PROP{s}, but for which only
compositions used in the definition of special elements (i.e.~the
`fractions') are allowed. The collection $\sfS \subset \Gamma(\Xi)$ would
then be the {\em free special\/} \PROP\ generated by $\Xi$.
}
\end{remark}

\begin{example}
\label{7}
{\rm
Let the boxes denote arbitrary special elements. Then the elements
\[
\frac{%
\frac{%
   \raisebox{.3em}{$\zeroatwo{} 
     \hskip .2em \zeroatwo{}$}}
{\twoatwo{}\hskip .2em\twoatwo {}} \hskip .3em
{
\unitlength=.8pt
\begin{picture}(50.00,30.00)(0.00,17.00)
\put(0.00,30.00){\line(0,-1){20.00}}
\put(50.00,30.00){\line(-1,0){50.00}}
\put(50.00,10.00){\line(0,1){20.00}}
\put(0.00,10.00){\line(1,0){50.00}}
\put(40.00,10.00){\line(0,-1){5.00}}
\put(30.00,10.00){\line(0,-1){5.00}}
\put(20.00,10.00){\line(0,-1){5.00}}
\put(10.00,10.00){\line(0,-1){5.00}}
\end{picture}}
}
{\twoazero {}\hskip .5em \twoazero {}\hskip .5em \twoazero {}
 \hskip .5em\twoazero {}}
\hskip 2mm \mbox { and } \hskip 2mm
\frac{\zeroatwo {}\hskip .5em {
\unitlength=.5pt
\begin{picture}(50.00,60.00)(0.00,10.00)
\put(30.00,10.00){\line(0,-1){5.00}}
\put(20.00,10.00){\line(0,-1){5.00}}
\put(50.00,50.00){\line(-1,0){20.00}}
\put(50.00,60.00){\line(0,-1){10.00}}
\put(30.00,60.00){\line(1,0){20.00}}
\put(30.00,50.00){\line(0,1){10.00}}
\put(0.00,60.00){\line(0,-1){10.00}}
\put(20.00,60.00){\line(-1,0){20.00}}
\put(20.00,50.00){\line(0,1){10.00}}
\put(0.00,50.00){\line(1,0){20.00}}
\put(30.00,30.00){\line(1,2){10.00}}
\put(20.00,30.00){\line(-1,2){10.00}}
\put(10.00,10.00){\line(0,1){20.00}}
\put(40.00,10.00){\line(-1,0){30.00}}
\put(40.00,30.00){\line(0,-1){20.00}}
\put(10.00,30.00){\line(1,0){30.00}}
\end{picture}}}{\twoazero {}\hskip 1em \twoazero {}}
\]
are also special, while the elements
\[
{
\unitlength=.5pt
\begin{picture}(110.00,60.00)(0.00,0.00)
\put(20.00,30.00){\line(0,-1){10.00}}
\put(10.00,30.00){\line(0,-1){10.00}}
\put(20.00,40.00){\line(0,-1){10.00}}
\put(10.00,40.00){\line(0,-1){10.00}}
\put(0.00,40.00){\line(0,1){20.00}}
\put(30.00,40.00){\line(-1,0){30.00}}
\put(30.00,60.00){\line(0,-1){20.00}}
\put(0.00,60.00){\line(1,0){30.00}}
\put(0.00,0.00){\line(0,1){20.00}}
\put(30.00,0.00){\line(-1,0){30.00}}
\put(30.00,20.00){\line(0,-1){20.00}}
\put(0.00,20.00){\line(1,0){30.00}}
\put(45.00,25.00){\makebox(0.00,0.00)[bl]{and}}
\put(30,0){
\put(90.00,40.00){\line(1,-2){10.00}}
\put(70.00,20.00){\line(1,0){40.00}}
\put(80.00,40.00){\line(0,-1){10.00}}
\put(70.00,40.00){\line(0,1){20.00}}
\put(110.00,40.00){\line(-1,0){40.00}}
\put(110.00,60.00){\line(0,-1){20.00}}
\put(70.00,60.00){\line(1,0){30.00}}
\put(110.00,20.00){\line(0,-1){20.00}}
\put(110.00,0.00){\line(-1,0){40.00}}
\put(70.00,0.00){\line(0,1){20.00}}
\put(80.00,30.00){\line(0,-1){10.00}}
\put(100.00,60.00){\line(1,0){10.00}}
\put(100.00,40.00){\line(-1,-2){10.00}}
}
\end{picture}}
\]
are not special. As an exercise, we recommend calculating the
genera of these composed elements in terms of the genera of individual
boxes. Other examples of special elements can be found in
Section~\ref{explicit}. 
}
\end{example}

The following lemma states that the path grading of special elements
from $\sfS(m,n)$ equals~$mn$.

\begin{lemma}
\label{strasna_vedra}
Let $m,n \geq 1$, let $X \in \sfS(m,n)$ be a monomial and let $1 \leq
i \leq m$, $1 \leq j \leq n$. Then there exists, in the graph $G_X$,
exactly one directed path connecting the $i$-th output with the $j$-th
input. In particular, $\pth(X) = mn$ for any $X \in \sfS(m,n)$.
\end{lemma}

\begin{proof}
The statement is certainly true for generators
$\xi^m_n$. Suppose we have proved it for some $\Rada A1l, \Rada B1k
\in \sfS$ and consider
\[
X : =\frac{A_1 \cdots A_l}{B_1 \cdots B_k}.
\]
There clearly exist unique $1 \leq s \leq l$ and $1 \leq t \leq k$ such
that the $i$-th output of $X$ is an output of $A_s$ and the $j$-th
input of $X$ is an input of $B_t$. 

It follows from the definition of $\sigma(k,l)$ that the $t$-th input
of $A_s$ is connected to the $s$-th output of $B_t$ and that the
$t$-th input of $A_s$ is the only input of $A_s$ which is connected to
some output of $B_t$. These considerations obviously imply that there
is, in $G_X$, a unique directed path connecting the $i$-th output with
the $j$-th input.%
\end{proof}

In the following lemma we give an upper bound for the genus of special
elements.

\begin{lemma}
\label{8}
Let $X \in \sfS(m,n)$ be a monomial. Then $\genus(X) \leq (m-1)(n-1)$.
\end{lemma}

\noindent
{\it Proof.} A straightforward induction on the `obvious' grading. If
$\grd(X) =1$, then $X$ is a generator and Lemma~\ref{8} holds
trivially.  
Each $X \in
\sfS(m,n)$ with $\grd(X) > 1$ can be decomposed as
\[
X = \frac {A_1 \cdots  A_v}{B_1 \cdots B_u},
\]
with some $1 \leq v \leq m$, $1 \le u \leq n$, $A_i \in \sfS(a_i,u)$,
$B_j \in \sfS(v,b_j)$, $a_i \geq 1$, $b_j \geq 1$, $1 \leq i \leq v$,
$1 \leq j \leq u$, $\sum_1^v a_i = m$, $\sum_1^u b_j = n$, such that
$\grd(A_i),\ \grd(B_j) < \grd(X)$. By Lemma~\ref{5} and the induction
assumption
\begin{eqnarray*}
\genus(X) &=&
(u-1)(v-1) + \textstyle\sum_1^v \genus(A_i) + \sum_1^u \genus(B_j)  
\\
\mbox { /by induction/ }
&\leq & (u-1)(v-1) +  \textstyle\sum_1^v (a_i -1)(u-1) + \sum_1^u (v-1)(b_j-1)
\\
&=& (u-1)(v-1) +  (m-v)(u-1) + (v-1)(n-u)
\\
&=&
(m-1)(n-1) - (m-v)(n-u) \leq (m-1)(n-1).\hfill  \qed
\end{eqnarray*}

\begin{remark}
\label{kasparek1}
{\rm\ Observe that the subspaces $\sfS(m,n) \subset \Gamma(\Xi)(m,n)$
are not $\Sigma_m$-$\Sigma_n$ invariant.  It easily follows from
Proposition~\ref{Zbraslavice_2003} and Lemma~\ref{strasna_vedra} that
the subspace $\sfS_0$ of $\sfS$ spanned by genus zero monomials
coincides with the free non-$\Sigma$ \hPROP\ $\uGamma_{\frac
12}(\uXi)$.  }
\end{remark}

\begin{theorem}
\label{jedeme_do_Jinolic}
Special elements form a friendly collection.
\end{theorem}

\begin{proof}
Condition (i) of Definition~\ref{s_tou_Elisabetou_jsem_to_prepisk},
with $C_{m,n} = (m-1)(n-1)$, follows from
Lemma~\ref{8}. Condition~(ii) follows from the fact, observed in
Remark~\ref{pocasi_nevychazi}, that $\sfS$ is the free special \PROP\
while (iii) is completely clear. In contrast, acyclicity~(iv) is a
very deep statement which we formulate as:%

\begin{proposition}
\label{11}
The vector spaces $\sfS(m,n)$ of special elements are $\pa_0$-acyclic
in positive degrees, for each $m,n \geq 1$.
\end{proposition}

Proposition~\ref{11} is proved in Section~\ref{ac}.
\end{proof}

\section{Explicit calculations}
\label{explicit}
\label{jeden_medvidek}

In this section we give a couple of formulas for the
perturbed differential (the formulas for the unperturbed differential
$\pa_0$ were given in Example~\ref{Orlik}). The first nontrivial
one expresses the compatibility axiom, the second two are~(\ref{A})
and~(\ref{BB}):  
\begin{eqnarray*}
\pa(\dvadva) &=& \pa_0(\dvadva) 
-
\frac{\jednadva \hskip .2em \jednadva}{\dvajedna \hskip .2em \dvajedna},
\\
\pa(\dvatri) &=& \pa_0(\dvatri) 
+
\frac{\jednadva \hskip .2em \jednadva}{\dvajedna \hskip .2em \dvadva}
-
\frac{\jednadva \hskip .2em \jednadva}{\dvadva \hskip .2em \dvajedna}
-
\frac{\bZbbZ \hskip .8em \jednatri}%
     {\dvajedna \hskip .2em \dvajedna \hskip .2em \dvajedna}
-
\frac{\jednatri \hskip .8em \ZbbZb}%
     {\dvajedna \hskip .2em \dvajedna \hskip .2em \dvajedna},
\\
\pa(\gen 32) &=& \pa_0(\tridva)
-
\frac{\jednadva \hskip .2em \dvadva}{\dvajedna \hskip .2em \dvajedna}
+
\frac{\dvadva \hskip .2em \jednadva}{\dvajedna \hskip .2em \dvajedna}
+
\frac{\jednadva \hskip .2em \jednadva \hskip .2em \jednadva}
     {\vZvvZ \hskip .8em \trijedna}
+
\frac{\jednadva \hskip .2em \jednadva \hskip .2em \jednadva}
     {\trijedna \hskip .8em \ZvvZv}.
\end{eqnarray*}

Let us pause a little and formulate the following conjecture.

\begin{conjecture}
\label{v_patek_prileti_Jim}
There exists a series of convex
$(m+n-3)$-dimensional polyhedra $B^m_n$ such that the differential
$\pa(\xi^m_n)$ is the sum of the codimension-one faces of these
polyhedra.
\end{conjecture}

\begin{figure}[t]
\begin{center}
{
\thinlines
\unitlength=.9pt
\begin{picture}(140.00,170.00)(0.00,-10.00)
\put(42.00,127.00){\makebox(0.00,0.00)[br]{$\frac{\gen 13 \hskip .5em \ZbbZb}%
                                    {\zdvihatko\gen 21\gen 21 \gen 21}$}}
\put(0.00,110.00){\makebox(0.00,0.00)[br]{$\frac{\bZbbZ \hskip .4em  \ZbbZb}%
                                         {\zdvihatko\gen 21\gen 21 \gen 21}$}}
\put(0.00,40.00){\makebox(0.00,0.00)[tr]{$\frac{\bZbbZ \hskip .4em\bZbbZ}%
                                       {\zdvihatko\gen 21 \gen 21 \gen 21}$}}
\put(0.00,70.00){\makebox(0.00,0.00)[r]{$\frac{\bZbbZ \hskip .5em \gen 13}%
                                         {\zdvihatko \gen 21\gen 21 \gen 21}$}}
\put(23.00,23.00){\makebox(0.00,0.00)[tr]{$\frac{\gen 12 \gen 12}%
                                        {\zdvihatko\gen 21\dvadva }$}}
\put(44.00,0.00){\makebox(0.00,0.00)[tr]{$\frac{\gen 12 \gen 12}%
                                    {\gen 21 \hskip .2em \dvojiteypsilon }$}}
\put(70.00,0.00){\makebox(0.00,0.00)[t]{$\dvabZbbZ$}}
\put(100.00,0.00){\makebox(0.00,0.00)[lt]{$\dvacarkabZbbZ$}}
\put(120.00,20.00){\makebox(0.00,0.00)[lt]{$\dvacarkatri$}}
\put(140.00,30.00){\makebox(0.00,0.00)[l]{$\dvacarkaZbbZb$}}
\put(140.00,73.00){\makebox(0.00,0.00)[l]{$\dvaZbbZb$}}
\put(140.00,110.00){\makebox(0.00,0.00)[l]{$\frac{\gen 12 \gen 12}%
                           {\hskip .1em\dvojiteypsilon \hskip .1em \gen 21}$}}
\put(107.00,127.00){\makebox(0.00,0.00)[lb]{$\frac{\gen 12 \gen 12}%
                                       {\zdvihatko\dvadva \gen 21}$}}
\put(70.00,145.00){\makebox(0.00,0.00)[b]{$\frac{\ZbbZb \hskip .2em\ZbbZb}%
                                     {\zdvihatko\gen 21 \gen 21 \gen 21}$}}
\put(10.00,100.00){\makebox(0.00,0.00){$\bullet$}}
\put(10.00,40.00){\makebox(0.00,0.00){$\bullet$}}
\put(40.00,10.00){\makebox(0.00,0.00){$\bullet$}}
\put(100.00,10.00){\makebox(0.00,0.00){$\bullet$}}
\put(130.00,40.00){\makebox(0.00,0.00){$\bullet$}}
\put(130.00,100.00){\makebox(0.00,0.00){$\bullet$}}
\put(70.00,130.00){\makebox(0.00,0.00){$\bullet$}}
\thicklines
\put(130.00,40.00){\line(-1,-1){30.00}}
\put(10.00,40.00){\line(1,-1){30.00}}
\put(40.00,10.00){\line(1,0){60.00}}
\put(10.00,100.00){\line(0,-1){60.00}}
\put(130.00,100.00){\line(0,-1){60.00}}
\put(70.00,130.00){\line(2,-1){60.00}}
\put(10.00,100.00){\line(2,1){60.00}}
\end{picture}}
\end{center}
\caption{\label{fig1}%
Heptagon $B^2_3$.
}
\end{figure}
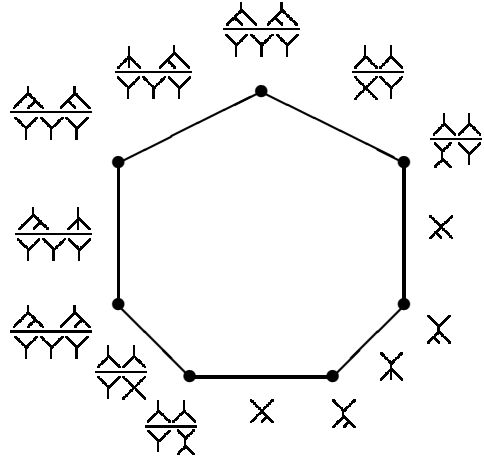

These polyhedra should generalize the case of $A_\infty$-algebras
discussed in Example~\ref{Orlik} in the sense that $B^1_n = B^n_1 =
K_n$ for $n \geq 2$. Clearly $B^2_2$ is the interval, while $B^2_3 =
B^3_2$ is the heptagon depicted in Figure~\ref{fig1}.
Before we proceed, we need to simplify our notation by an almost
obvious `linear extension' of $(k,l)$-fractions.

\begin{notation}
{\rm\
Let $k,l,s,t \geq 1$, $\Rada {A^s}1l \in \sfS(*,k)$ and $\Rada {B^t}1k \in
\sfS(l,*)$.
Then define
\[
\frac{\textstyle\sum_s A^s_1 \cdots A^s_l}%
     {\textstyle\sum_t B^t_1 \cdots B^t_k} :=
\sum_{s,t} \frac{A^s_1 \cdots A^s_l}{B^t_1 \cdots B^t_k}~.
\]
}
\end{notation}

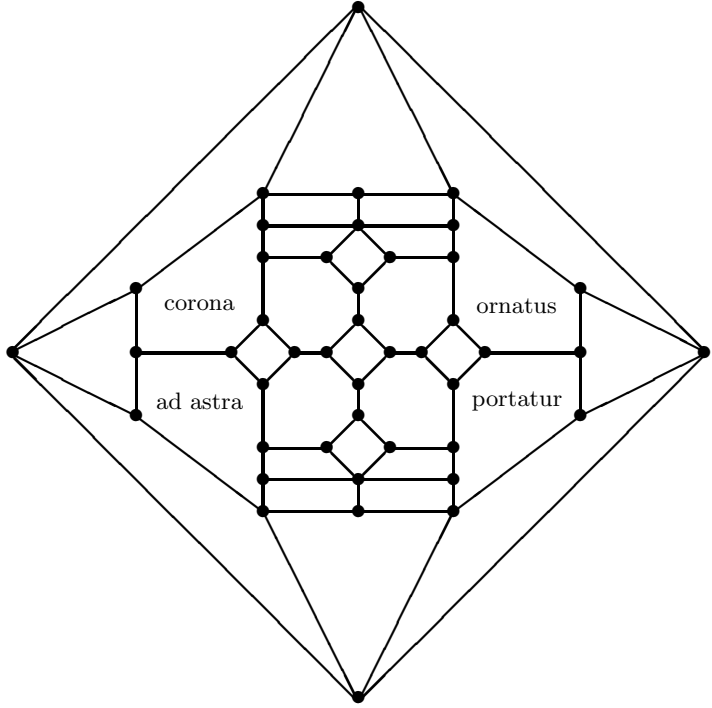
\begin{figure}[t]
\begin{center}
{
\unitlength=0.6pt
\thicklines
\begin{picture}(440.00,435.00)(0.00,0.00)
\put(320.00,190.00){\makebox(0.00,0.00){\scriptsize portatur}}
\put(120.00,190.00){\makebox(0.00,0.00){\scriptsize ad astra}}
\put(320.00,250.00){\makebox(0.00,0.00){\scriptsize ornatus}}
\put(120.00,250.00){\makebox(0.00,0.00){\scriptsize corona}}
\put(280.00,120.00){\makebox(0.00,0.00){$\bullet$}}
\put(160.00,120.00){\makebox(0.00,0.00){$\bullet$}}
\put(220.00,140.00){\makebox(0.00,0.00){$\bullet$}}
\put(280.00,140.00){\makebox(0.00,0.00){$\bullet$}}
\put(160.00,140.00){\makebox(0.00,0.00){$\bullet$}}
\put(280.00,160.00){\makebox(0.00,0.00){$\bullet$}}
\put(240.00,160.00){\makebox(0.00,0.00){$\bullet$}}
\put(200.00,160.00){\makebox(0.00,0.00){$\bullet$}}
\put(160.00,160.00){\makebox(0.00,0.00){$\bullet$}}
\put(360.00,180.00){\makebox(0.00,0.00){$\bullet$}}
\put(280.00,200.00){\makebox(0.00,0.00){$\bullet$}}
\put(160.00,200.00){\makebox(0.00,0.00){$\bullet$}}
\put(220.00,200.00){\makebox(0.00,0.00){$\bullet$}}
\put(220.00,180.00){\makebox(0.00,0.00){$\bullet$}}
\put(80.00,180.00){\makebox(0.00,0.00){$\bullet$}}
\put(360.00,220.00){\makebox(0.00,0.00){$\bullet$}}
\put(300.00,220.00){\makebox(0.00,0.00){$\bullet$}}
\put(260.00,220.00){\makebox(0.00,0.00){$\bullet$}}
\put(240.00,220.00){\makebox(0.00,0.00){$\bullet$}}
\put(200.00,220.00){\makebox(0.00,0.00){$\bullet$}}
\put(180.00,220.00){\makebox(0.00,0.00){$\bullet$}}
\put(140.00,220.00){\makebox(0.00,0.00){$\bullet$}}
\put(80.00,220.00){\makebox(0.00,0.00){$\bullet$}}
\put(360.00,260.00){\makebox(0.00,0.00){$\bullet$}}
\put(80.00,260.00){\makebox(0.00,0.00){$\bullet$}}
\put(280.00,240.00){\makebox(0.00,0.00){$\bullet$}}
\put(220.00,260.00){\makebox(0.00,0.00){$\bullet$}}
\put(220.00,240.00){\makebox(0.00,0.00){$\bullet$}}
\put(160.00,240.00){\makebox(0.00,0.00){$\bullet$}}
\put(280.00,280.00){\makebox(0.00,0.00){$\bullet$}}
\put(240.00,280.00){\makebox(0.00,0.00){$\bullet$}}
\put(200.00,280.00){\makebox(0.00,0.00){$\bullet$}}
\put(160.00,280.00){\makebox(0.00,0.00){$\bullet$}}
\put(280.00,300.00){\makebox(0.00,0.00){$\bullet$}}
\put(220.00,300.00){\makebox(0.00,0.00){$\bullet$}}
\put(160.00,300.00){\makebox(0.00,0.00){$\bullet$}}
\put(280.00,320.00){\makebox(0.00,0.00){$\bullet$}}
\put(160.00,320.00){\makebox(0.00,0.00){$\bullet$}}
\put(220.00,120.00){\makebox(0.00,0.00){$\bullet$}}
\put(220.00,320.00){\makebox(0.00,0.00){$\bullet$}}
\put(438.00,220.00){\makebox(0.00,0.00){$\bullet$}}
\put(2.00,220.00){\makebox(0.00,0.00){$\bullet$}}
\put(220.00,2.00){\makebox(0.00,0.00){$\bullet$}}
\put(220.00,438.00){\makebox(0.00,0.00){$\bullet$}}
\put(220.00,0.00){\line(1,2){60.00}}
\put(160.00,120.00){\line(1,-2){60.00}}
\put(220.00,0.00){\line(-1,1){220.00}}
\put(440.00,220.00){\line(-1,-1){220.00}}
\put(80.00,180.00){\line(4,-3){80.00}}
\put(440.00,220.00){\line(-2,-1){80.00}}
\put(360.00,180.00){\line(-4,-3){80.00}}
\put(360.00,220.00){\line(0,-1){40.00}}
\put(220.00,440.00){\line(1,-1){220.00}}
\put(360.00,260.00){\line(2,-1){80.00}}
\put(360.00,260.00){\line(-4,3){80.00}}
\put(360.00,220.00){\line(0,1){40.00}}
\put(300.00,220.00){\line(1,0){60.00}}
\put(220.00,440.00){\line(1,-2){60.00}}
\put(80.00,220.00){\line(1,0){30.00}}
\put(80.00,180.00){\line(0,1){10.00}}
\put(0.00,220.00){\line(2,-1){80.00}}
\put(0.00,220.00){\line(2,1){80.00}}
\put(80.00,260.00){\line(0,-1){70.00}}
\put(80.00,260.00){\line(4,3){80.00}}
\put(220.00,440.00){\line(-1,-1){220.00}}
\put(160.00,320.00){\line(1,2){60.00}}
\put(110.00,220.00){\line(1,0){30.00}}
\put(240.00,280.00){\line(1,0){40.00}}
\put(160.00,280.00){\line(1,0){40.00}}
\put(240.00,160.00){\line(1,0){40.00}}
\put(160.00,160.00){\line(1,0){40.00}}
\put(160.00,300.00){\line(1,0){120.00}}
\put(160.00,140.00){\line(1,0){120.00}}
\put(160.00,120.00){\line(0,1){80.00}}
\put(280.00,120.00){\line(-1,0){120.00}}
\put(280.00,200.00){\line(0,-1){80.00}}
\put(280.00,320.00){\line(0,-1){80.00}}
\put(160.00,320.00){\line(1,0){120.00}}
\put(160.00,240.00){\line(0,1){80.00}}
\put(160.00,200.00){\line(1,1){20.00}}
\put(140.00,220.00){\line(1,-1){20.00}}
\put(160.00,240.00){\line(-1,-1){20.00}}
\put(180.00,220.00){\line(-1,1){20.00}}
\put(200.00,220.00){\line(-1,0){20.00}}
\put(240.00,160.00){\line(-1,1){20.00}}
\put(220.00,140.00){\line(1,1){20.00}}
\put(200.00,160.00){\line(1,-1){20.00}}
\put(220.00,180.00){\line(-1,-1){20.00}}
\put(220.00,200.00){\line(0,-1){20.00}}
\put(220.00,140.00){\line(0,-1){20.00}}
\put(220.00,320.00){\line(0,-1){20.00}}
\put(280.00,240.00){\line(-1,-1){20.00}}
\put(300.00,220.00){\line(-1,1){20.00}}
\put(280.00,200.00){\line(1,1){20.00}}
\put(260.00,220.00){\line(1,-1){20.00}}
\put(240.00,220.00){\line(1,0){20.00}}
\put(200.00,280.00){\line(1,-1){20.00}}
\put(220.00,300.00){\line(-1,-1){20.00}}
\put(240.00,280.00){\line(-1,1){20.00}}
\put(220.00,260.00){\line(1,1){20.00}}
\put(220.00,260.00){\line(0,-1){20.00}}
\put(220.00,240.00){\line(0,1){20.00}}
\put(220.00,240.00){\line(0,1){20.00}}
\put(220.00,240.00){\line(0,1){20.00}}
\put(220.00,240.00){\line(0,1){20.00}}
\put(220.00,200.00){\line(-1,1){20.00}}
\put(240.00,220.00){\line(-1,-1){20.00}}
\put(220.00,240.00){\line(1,-1){20.00}}
\put(200.00,220.00){\line(1,1){20.00}}
\end{picture}}
\end{center}
\caption{\label{fig2}%
The plane projection of 3-dimensional polyhedron $B^3_3$ from one
of its square faces. Polyhedron $B^3_3$
has 30 2-dimensional 
faces (8 heptagons and 22 squares), 72 edges and 44 vertices.}
\end{figure}

\noindent 
For example, with this notation the formula for $\pa(\dvatri)$ can be
simplified to
\begin{eqnarray*}
\pa(\dvatri) &=& \pa_0(\dvatri) 
+
\frac{\jednadva \hskip .2em \jednadva}
{\dvajedna \hskip .2em \dvadva - \dvadva \hskip .2em \dvajedna}
-
\frac{\bZbbZ \hskip .2em \jednatri + \jednatri \hskip .2em \ZbbZb}%
     {\dvajedna \hskip .2em \dvajedna \hskip .2em \dvajedna}
\\
&=&
 \pa_0(\dvatri) 
+
\frac{\Delta(\jednadva)}
{\dvajedna \hskip .2em \dvadva - \dvadva \hskip .2em \dvajedna}
-
\frac{\Delta (\gen 13)}{\dvajedna \hskip .2em \dvajedna \hskip .2em
\dvajedna}
\hskip 2pt,
\end{eqnarray*}
where $\Delta$ is the Saneblidze-Umble diagonal in the
associahedron~\cite{saneblidze-umble:preprint2}. 
\eject
The next term is
\begin{eqnarray*}
\pa(\tritri) &=& \pa_0(\tritri)
-
\frac{\gen12 \hskip .2em \gen22}{\gen22 \hskip .2em \gen21}
-
\frac{\gen22 \hskip .2em \gen12}{\gen21 \hskip .2em \gen22}
+
\frac{\gen22 \hskip .2em \gen12}{\gen22 \hskip .2em \gen21}
+
\frac{\gen12 \hskip .2em \gen22}{\gen21 \hskip .2em \gen22}
\\
&& \hskip -.7cm
+
\frac{\gen12 \hskip .2em \gen12 \hskip .2em \gen12}%
     {\gen32 \hskip .3em \ZvvZv -
     \vZvvZ \hskip .3em \gen32  -
     \vZvvZdva \hskip .3em \gen31 -
     \gen31 \hskip .3em \ZvvZvdva}
-
\frac{\bZbbZ \hskip .3em \gen23 
      -\gen23 \hskip .3em \ZbbZb - 
     \dvabZbbZ \hskip .3em \gen13
     - \gen13 \hskip .3em \dvaZbbZb}%
     {\gen21 \hskip .2em \gen21 \hskip .2em \gen21}
\\
&&\hskip -.7cm
-
\frac{\gen12 \hskip .2em \gen12 \hskip .2em \gen12}%
     {\gen31 \hskip .3em \frac{\gen22 \hskip .2em \gen 12}%
     {\rule{0pt}{9pt}\gen 21 \hskip .2em \gen21}
      + \frac{\gen12 \hskip .2em \gen 22}%
     {\rule{0pt}{9pt}\gen 21 \hskip .2em \gen21} \hskip .3em \gen31}
+
\frac {\gen13 \hskip .3em \frac{\gen12 \hskip .2em \gen 12}%
     {\rule{0pt}{9pt}\gen 22 \hskip .2em \gen21}
      + \frac{\gen12 \hskip .2em \gen 12}%
     {\rule{0pt}{9pt}\gen 21 \hskip .2em \gen22} \hskip .3em \gen13}%
     {\gen21 \hskip .2em \gen21 \hskip .2em \gen21}
-
\frac{\bZbbZ  \hskip .2em \bZbbZ \hskip .2em \gen13
      + \bZbbZ \hskip .2em \gen13 \hskip .2em \ZbbZb
      + \gen13 \hskip .2em \ZbbZb \hskip .2em \ZbbZb 
     }
     {\vZvvZ  \hskip .2em \vZvvZ \hskip .2em \gen31
     + \vZvvZ  \hskip .2em  \gen31 \hskip .2em \ZvvZv
     + \gen31 \hskip .2em \ZvvZv   \hskip .2em \ZvvZv
     }\hskip 1mm.
\end{eqnarray*}
Observe that the last term of the above equation is
\[
-\frac{\Delta^{(3)}(\gen 13)}{\Delta^{(3)}(\gen 31)},
\]
where $\Delta^{(3)}(-) := (\Delta \ot \id)\Delta(-)$ 
denotes the iteration of the
Saneblidze-Umble diagonal which is coassociative on $\gen 13$ and
$\gen 31$ (see~\cite{markl-shnider:dg}).
The corresponding 3-dimensional polyhedron $B^3_3$ is shown in
Figure~\ref{fig2}.

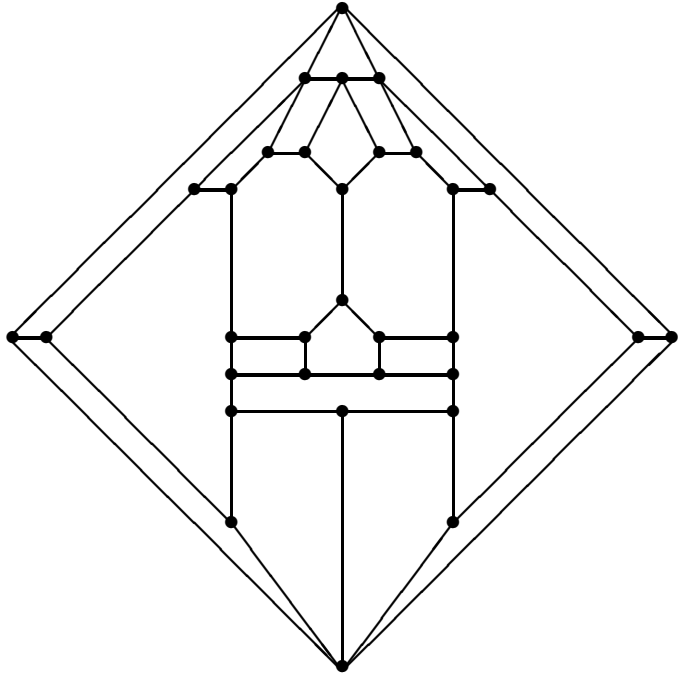
\begin{figure}[t]
\begin{center}
{
\unitlength=.7pt
\begin{picture}(360.00,360.00)(0.00,0.00)
\thicklines
\put(180.00,2.00){\makebox(0.00,0.00){$\bullet$}}
\put(240.00,80.00){\makebox(0.00,0.00){$\bullet$}}
\put(120.00,80.00){\makebox(0.00,0.00){$\bullet$}}
\put(240.00,140.00){\makebox(0.00,0.00){$\bullet$}}
\put(120.00,140.00){\makebox(0.00,0.00){$\bullet$}}
\put(240.00,160.00){\makebox(0.00,0.00){$\bullet$}}
\put(200.00,160.00){\makebox(0.00,0.00){$\bullet$}}
\put(160.00,160.00){\makebox(0.00,0.00){$\bullet$}}
\put(120.00,160.00){\makebox(0.00,0.00){$\bullet$}}
\put(358.00,180.00){\makebox(0.00,0.00){$\bullet$}}
\put(340.00,180.00){\makebox(0.00,0.00){$\bullet$}}
\put(240.00,180.00){\makebox(0.00,0.00){$\bullet$}}
\put(200.00,180.00){\makebox(0.00,0.00){$\bullet$}}
\put(160.00,180.00){\makebox(0.00,0.00){$\bullet$}}
\put(120.00,180.00){\makebox(0.00,0.00){$\bullet$}}
\put(20.00,180.00){\makebox(0.00,0.00){$\bullet$}}
\put(2.00,180.00){\makebox(0.00,0.00){$\bullet$}}
\put(180.00,200.00){\makebox(0.00,0.00){$\bullet$}}
\put(260.00,260.00){\makebox(0.00,0.00){$\bullet$}}
\put(240.00,260.00){\makebox(0.00,0.00){$\bullet$}}
\put(180.00,260.00){\makebox(0.00,0.00){$\bullet$}}
\put(120.00,260.00){\makebox(0.00,0.00){$\bullet$}}
\put(100.00,260.00){\makebox(0.00,0.00){$\bullet$}}
\put(220.00,280.00){\makebox(0.00,0.00){$\bullet$}}
\put(200.00,280.00){\makebox(0.00,0.00){$\bullet$}}
\put(160.00,280.00){\makebox(0.00,0.00){$\bullet$}}
\put(140.00,280.00){\makebox(0.00,0.00){$\bullet$}}
\put(180.00,320.00){\makebox(0.00,0.00){$\bullet$}}
\put(200.00,320.00){\makebox(0.00,0.00){$\bullet$}}
\put(160.00,320.00){\makebox(0.00,0.00){$\bullet$}}
\put(180.00,358.00){\makebox(0.00,0.00){$\bullet$}}
\put(180.00,140.00){\makebox(0.00,0.00){$\bullet$}}
\put(180.00,0.00){\line(0,1){140}}
\put(180.00,0.00){\line(3,4){60.00}}
\put(120.00,80.00){\line(3,-4){60.00}}
\put(180.00,0.00){\line(1,1){180.00}}
\put(0.00,180.00){\line(1,-1){180.00}}
\put(120.00,80.00){\line(0,1){60.00}}
\put(20.00,180.00){\line(1,-1){100.00}}
\put(240.00,80.00){\line(0,1){60.00}}
\put(340.00,180.00){\line(-1,-1){100.00}}
\put(120.00,260.00){\line(-1,0){20.00}}
\put(20.00,180.00){\line(1,1){140.00}}
\put(0.00,180.00){\line(1,0){20.00}}
\put(180.00,360.00){\line(-1,-1){180.00}}
\put(240.00,260.00){\line(1,0){20.00}}
\put(360.00,180.00){\line(-1,0){20.00}}
\put(180.00,360.00){\line(1,-1){180.00}}
\put(200.00,320.00){\line(1,-1){140.00}}
\put(180.00,360.00){\line(1,-2){20.00}}
\put(160.00,320.00){\line(1,2){20.00}}
\put(160.00,320.00){\line(-1,-2){20.00}}
\put(170.00,320.00){\line(-1,0){10.00}}
\put(200.00,320.00){\line(1,-2){20.00}}
\put(170.00,320.00){\line(1,0){30.00}}
\put(180.00,320.00){\line(1,-2){20.00}}
\put(160.00,280.00){\line(1,2){20.00}}
\put(120.00,260.00){\line(0,-1){80.00}}
\put(140.00,280.00){\line(-1,-1){20.00}}
\put(160.00,280.00){\line(-1,0){20.00}}
\put(240.00,260.00){\line(0,-1){80.00}}
\put(220.00,280.00){\line(1,-1){20.00}}
\put(200.00,280.00){\line(1,0){20.00}}
\put(180.00,260.00){\line(-1,1){20.00}}
\put(180.00,260.00){\line(1,1){20.00}}
\put(180.00,220.00){\line(0,1){40.00}}
\put(180.00,200.00){\line(0,1){20.00}}
\put(240.00,140.00){\line(0,1){20.00}}
\put(120.00,140.00){\line(1,0){120.00}}
\put(120.00,160.00){\line(0,-1){20.00}}
\put(240.00,180.00){\line(-1,0){40.00}}
\put(240.00,160.00){\line(0,1){20.00}}
\put(200.00,160.00){\line(1,0){40.00}}
\put(120.00,160.00){\line(1,0){40.00}}
\put(120.00,180.00){\line(0,-1){20.00}}
\put(160.00,180.00){\line(-1,0){40.00}}
\put(180.00,200.00){\line(-1,-1){20.00}}
\put(200.00,180.00){\line(-1,1){20.00}}
\put(200.00,160.00){\line(0,1){20.00}}
\put(160.00,160.00){\line(1,0){40.00}}
\put(160.00,180.00){\line(0,-1){20.00}}
\end{picture}}
\end{center}
\caption{\label{fig3}%
The plane projection of 3-dimensional polyhedron $B^2_4$. It has 32
vertices, 51 edges and 21 two-dimensional faces (3 pentagons, 
5 heptagons and 13 squares).}
\end{figure}

\eject
The relation with the Saneblidze-Umble diagonal $\Delta$ is even more manifest
in the formula

\begin{eqnarray*}
\pa(\dvactyri) &=& \pa_0(\dvactyri) 
+
\frac{\gen 12 \gen 12}{\gen 23 \gen 21 + \gen 22 \gen 22 + \gen 21
\gen 23}
+ 
\frac{\bZbbZ \gen 13 + \gen 13 \ZbbZb}{\gen 22 \gen 21 \gen 21 - \gen
21 \gen 22 \gen 21 + \gen 21 \gen 21 \gen 22}
\\
&&
+
\frac{
      \bZbbbZ \bZbbZb + \bZbbZb \ZbbbZb + \bZbbbZ \ZbbbZb 
       -\bZbZbbZ \gen 14 - \bbZbbZ \ZbbZbb - \gen 14 \ZbbZbZb}
      {\gen 21 \gen 21 \gen 21 \gen 21}
\\
&=& \pa_0(\dvactyri) 
+
\frac{\Delta(\gen 12)}{\gen 23 \gen 21 + \gen 22 \gen 22 + \gen 21
\gen 23}
+ 
\frac{\Delta(\gen 13)}{\gen 22 \gen 21 \gen 21 - \gen
21 \gen 22 \gen 21 + \gen 21 \gen 21 \gen 22}
+
\frac{\Delta(\gen 14)}{\gen 21 \gen 21 \gen 21 \gen 21}.
\end{eqnarray*}
The corresponding $B^2_4$ is shown in Figure~\ref{fig3}.

\section{Calculus of special elements}
\label{calc}
\label{kam_nas_vystehuji?}

This section provides a preparatory material for the
proof of the $\pa_0$-acyclicity of the space
$\sfS(m,n)$ given in Section~\ref{ac}. 
As in the proof of Lemma~\ref{8}, each monomial $X \in
\sfS(m,n)$ is represented as
\begin{equation}
\label{eq2}
X = \frac{A_1 \cdots A_v}{B_1 \cdots B_u},
\end{equation}
for $1 \leq v \leq m$, $1 \leq u  \leq n$, $A_i \in \sfS(a_i,u)$, $B_j \in
\sfS(v,b_j)$, with $\sum_1^v a_i =m$ and $\sum_1^u b_j =n$.
Very crucially, representation~(\ref{eq2}) is \underline{not}
\underline{uni}q\underline{ue}, as illustrated in the following example.

\begin{example}
{\rm\
\label{12}
It is easy to verify that
\begin{equation}
\label{obscese}
\frac{\jednadva \hskip 1.3em \jednadva}
{\frac{\jednadva \hskip .2em \jednadva}{\rule{0pt}{9pt}\dvajedna \hskip .2em
\dvajedna} \hskip .2em \dvajedna} 
=
\frac{\ZbbZb \hskip .5em  \ZbbZb}
{\dvajedna \hskip .2em \dvajedna\hskip .2em \dvajedna } \hskip 2pt.
\end{equation}
Therefore, the element $X \in \sfS(2,3)$ above can be
either  represented as
\[
X = \frac{A_1 A_2}{B_1 B_2},
\]
with $A_1 = A_2 = \jednadva \in \sfS(1,2)$, 
\[
B_1 = \frac{\jednadva \hskip .2em \jednadva}{\dvajedna \hskip .2em
\dvajedna} \in \sfS(2,2)
\]
and $B_2 = \dvajedna \in \sfS(2,1)$, or as
\[
X = \frac{A'_1 \hskip 1.2em A'_2}{B'_1 B'_2 B'_3},
\]
with
\[
A'_1 = A'_2 = \ZbbZb \in \sfS(1,3)
\]
and $B'_1 = B'_2  = B'_3 = \dvajedna \in \sfS(2,1)$.
Of a bit different nature is the relation
\begin{equation}
\label{caj}
\frac{
\frac{\jednadva \hskip .2em \jednadva}
{\hskip .1em\rule{0pt}{10pt}\twotwo \hskip .3em \dvajedna}
\hskip .2em \ZbbZb
}
{
\dvajedna \hskip .4em \dvajedna \hskip .4em \dvajedna
}
=
\frac{\jednadva \hskip .4em \jednadva \hskip .4em \jednadva}
{
\frac{\twotwo \hskip .2em \jednadva}
     {\hskip .1em \rule{0pt}{9pt}\dvajedna \hskip .3em \dvajedna}
\hskip .2em \ZvvZv}
\end{equation}
or a similar one
\begin{equation}
\label{caj_zeleny}
\frac{
\frac{\jednadva \hskip .2em \jednadva}%
{\hskip -.1em\rule{0pt}{10pt}\dvajedna \hskip .3em
\twotwo}
\hskip .2em
\bZbbZ
}
{
\dvajedna \hskip .4em \dvajedna \hskip .4em \dvajedna
}
=
\frac{\jednadva \hskip .4em \jednadva \hskip .4em \jednadva}
{\ZvvZv \hskip .2em
\frac{\hskip -.1em \twotwo \hskip .3em \jednadva}
     {\rule{0pt}{9pt}\dvajedna \hskip .2em \dvajedna}
}\hskip 2pt,
\end{equation}
where \raisebox{-2pt}{$\twotwo$} is an arbitrary element of $\sfS(2,2)$.
}
\end{example}
It follows from the above remarks that
\begin{equation}
\label{uz_je_hodne_pozde}
\sfS(m,n) 
=  \span(\xi^m_n) \oplus  
\bigoplus_M \sfS(a_1,u) \ot \cdots \ot \sfS(a_v,u)
\ot
\sfS(v,b_1) \ot \cdots \ot \sfS(v,b_u)/R(m,n)
\end{equation}
where
\[
M = \{
1\leq v \leq m,\ 1 \leq u  \leq n,\ 
(v,u) \not= (1,n),(m,1) \mbox { and }
\textstyle\sum_1^v a_i =m,\ \sum_1^u b_j =n
\}
\]
and $R(m,n)$ accounts for the non-uniqueness of presentation~(\ref{eq2}). 
Observe that if $R(m,n)$ were trivial, then the
$\pa_0$-acyclicity of $\sfS(m,m)$ would follow immediately from the K\"unneth
formula and induction.

\begin{example}
{\rm\
We have
\begin{eqnarray*}
\sfS(2,2)& \cong & \span(\xi^2_2) \oplus [\S(2,1) \otimes \S(1,2)] \oplus
[\S(1,2) \ot \S(1,2) \ot \S(2,1) \ot \S(2,1)]
\\
&\cong &
\span(\dvadva) \oplus \span(\dvojiteypsilon) \oplus 
\span\left(
    \frac{\jednadva \hskip .2em \jednadva}{\dvajedna \hskip .2em\dvajedna}
    \right),
\end{eqnarray*}
the relations $R(2,2)$ are trivial. On the other hand, the left-hand
side of~(\ref{caj}) represents an element of $\S(3,3)$ by
\[
\frac{\jednadva \hskip .2em \jednadva}{\rule{0pt}{10pt}\twotwo \hskip .2em
\dvajedna} \ot  \ZbbZb \ot \dvajedna \ot \dvajedna \ot \dvajedna 
\in \S(2,3) \ot \S(1,3) \ot \S(2,1)\ot \S(2,1)\ot \S(2,1),
\]
while the right-hand side of~(\ref{caj}) represents the same element by
\[
\jednadva \ot \jednadva \ot \jednadva 
\ot
\frac{\twotwo \hskip .2em \jednadva}
     {\rule{0pt}{9pt}\dvajedna \hskip .2em \dvajedna}
\ot \ZvvZv 
\in \S(1,2) \ot\S(1,2) \ot\S(1,2) \ot \S(3,2) \ot \S(3,1).
\]
Therefore $R(3,3)$ must contain a relation that identifies these two elements.
}
\end{example}

Let us describe the space of relations $R$.  Suppose that $s,t \geq
1$, $\Rada c1s,\Rada d1t \geq 1$ are natural numbers and let \cd\
denote the array $(\Rada c1s;\Rada d1t)$. A crucial r\^ole in the
following definition will be played by a matrix
\begin{equation}
\label{potim_se}
C = 
\left(C_{ij} \right)_{\hskip -.3em \doubless{1 \leq i \leq t}{1 \leq j \leq s}}
\end{equation}
with entries $C_{ij} \in \sfS(d_i,c_j)$. Finally, let
\[
X = \frac{A_1 \cdots A_v}{B_1 \cdots B_u}
\]
be a monomial as in~(\ref{eq2}).

\begin{definition}
Element $X$ is called {\em \cd-up-reducible\/} if $u = c_1 + \cdots +
c_s$, $v=t$ and
\begin{equation}
\label{U1}
A_i = \frac {A_{i1} \cdots A_{id_i}}{C_{i1} \cdots C_{is}}
\end{equation}
for some $A_{ib} \in \sfS(*,s)$, $1 \leq i \leq t$, $1 \leq b \leq
d_i$, where $C_{ij}$ are entries of a matrix as in~(\ref{potim_se}).
Dually, $X$ is called {\em \cd-down-reducible\/} if $v = d_1 + \cdots
+ d_t$, $u=s$ and
\begin{equation}
\label{U2}
B_j = \frac {C_{1j} \cdots C_{tj}}{B_{1j} \cdots B_{c_jj}}
\end{equation}
for some $B_{aj} \in \sfS(t,*)$, $1 \leq j \leq s$, $1 \leq a \leq
c_j$, again with $C = (C_{ij})$ as in~(\ref{potim_se}).  We
denote by $\Up\cd$ (resp.~\Dw\cd) the subspace spanned by all \cd-
(resp.~down) -reducible monomials.
\end{definition}

\begin{proposition}
\label{opravdu_pozde}
The spaces $\Up\cd$ and $\Dw\cd$ are isomorphic. The isomorphism
is given by the identification of the up-reducible element
\[
\frac{
\hbox{\large$\frac{A_{11} \cdots A_{1d_1}}{C_{11} \cdots C_{1s}}
\hskip .8em
\frac{A_{21} \cdots A_{2d_2}}{C_{21} \cdots C_{2s}}\hskip .8em
\cdots \hskip .8em
\frac{A_{t1} \cdots A_{td_t}}{C_{t1} \cdots C_{ts}}$}
}
{
B_{11} \cdots B_{c_11} \hskip .4em B_{12} 
\cdots B_{c_22} \hskip .2em   \cdots \hskip .2em B_{1s}\cdots B_{c_ss}
}
\]
with the down-reducible element
\[
\frac{
A_{11} \cdots A_{1d_1} \hskip .4em A_{21} 
\cdots A_{2d_2} \hskip .2em   \cdots \hskip .2em A_{t1}\cdots A_{td_t}
}
{
\hbox{\large$\frac{C_{11} \cdots C_{t1}}{B_{11} \cdots B_{c_11}}
\hskip .8em
\frac{C_{12} \cdots C_{t2}}{B_{12} \cdots B_{c_22}}\hskip .8em
\cdots \hskip .8em
\frac{C_{1s} \cdots C_{ts}}{B_{1s} \cdots B_{c_ss}}$}
}\hskip 2pt .
\]
Relations $R$ in~(\ref{uz_je_hodne_pozde}) are generated by the
above identifications. 
\end{proposition}

\begin{proof}
The proof follows from analyzing the underlying graphs.%
\end{proof}

We call the relations described in
Proposition~\ref{opravdu_pozde} the {\em \cd-relations\/}.
These relations are clearly compatible with the differential $\pa_0$
and do not change the genus. They are trivial if $d_j = c_i =1$, 
for all $i,j$. 

\begin{example}
{\rm\ Equation~(\ref{obscese}) of Example~\ref{12} is an equality of
two $(2,1;1,1)$-reducible elements with $A_{11} = A_{21} = \jednadva$,
$B_{11} = B_{21} = B_{12} = \dvajedna$ \hskip .3em and the
matrix~(\ref{potim_se}) given by
\[
C =
\left(
\begin{array}{cc}
\gen12 & \id
\\
\gen12 & \id
\end{array}
\right).
\]
Equation~(\ref{caj}) of Example~\ref{12} is an equality between two
$(2,1;2,1)$-reducible elements with $A_{11} = A_{12} = A_{21} =
\jednadva$, $B_{11} = B_{21} = B_{12} = \dvajedna$ and
\[
C =
\left(
\begin{array}{cc}
\twotwo & \gen21
\\
\gen12 & \id
\end{array}
\right).
\]
We leave it as an exercise to interpret also~(\ref{caj_zeleny}) in
terms of \cd-relations.
}
\end{example}

\begin{example}
\label{kourim_vodni_dymku}
{\rm\
Let us write presentation~(\ref{uz_je_hodne_pozde}) for $\S(2,3)$. Of
course, $\S(2,3) = \S_0(2,3) \op \S_1(2,3) \op \S_2(2,3)$, where the
subscript denotes the genus. Then $\S_0(2,3)$ is represented as the
quotient of
\begin{eqnarray*}
&& \hskip -1cm \Span(\xi^2_3) \op  [\S_0(2,2) \!\ot\! \S(1,1) \!\ot\! \S(1,2)]
\op [\S(2,1) \!\ot\! \S(1,3)] \op [\S_0(2,2) \!\ot\! \S(1,2)
\!\ot\! \S(1,1)] \cong
\\
&&\cong
\span(\dvatri) 
\op 
\span\left(\frac{\hskip 3pt\raisebox{-2pt}{\twotwo}}%
                {\id \hskip .3em\gen12}\right) 
\op  
\span\left(\frac{\hskip -8pt\raisebox{-2pt}{\gen 21}}%
{\rule{0pt}{14pt}\onethree}\right)
\op 
\span\left(\frac{\hskip 3pt\raisebox{-2pt}{\twotwo}}%
{\gen 12 \hskip .3em \id}\right),
\end{eqnarray*}
where \raisebox{-2pt}{$\onethree$} 
is an arbitrary element of $\S_0(1,3)$ and \raisebox{-2pt}{$\twotwo$} is
an element of $\S_0(2,2)$, modulo relations
$R(2,3)$ that identify the up-$(2;1)$-reducible element
\[
\frac{\hskip 3pt\raisebox{-2pt}{\dvojiteypsilon}}%
                {\id \hskip .3em\gen12}
\in
\span\left(\frac{\hskip 3pt\raisebox{-2pt}{\twotwo}}%
                {\id \hskip .3em\gen12}\right)
\] 
with the down-$(2;1)$-reducible element
\[
\frac{\hskip -8pt\raisebox{-2pt}{\gen 21}}{\bZbbZ}
\in
\span\left(\frac{\hskip -8pt\raisebox{-2pt}{\gen 21}}%
{\rule{0pt}{14pt}\onethree}\right)
\]
and identify the up-$(2;1)$-reducible element
\[
\frac{\hskip 3pt\raisebox{-2pt}{\dvojiteypsilon}}%
                {\gen 12 \hskip .3em\id}
\in
\span\left(\frac{\hskip 3pt\raisebox{-2pt}{\twotwo}}%
                {\gen 12 \hskip .3em\id}\right) 
\]
with the down-$(2;1)$-reducible element
\[
\frac{\hskip -8pt\raisebox{-2pt}{\gen 21}}{\ZbbZb}
\in
\span\left(\frac{\hskip -8pt\raisebox{-2pt}{\gen 21}}%
{\rule{0pt}{14pt}\onethree}\right).
\]
With the obvious similar notation, $\S_1(2,3)$ is the quotient of
\[
\span\left(
\frac{\gen 12 \hskip .2em \gen 12}{\twotwo \hskip .2em \gen 21}
\right) 
\op
\span\left(\frac{\frac{\gen 12 \gen 12}%
{\rule{0pt}{9pt}\gen 21 \gen 21}}%
{\hskip 5pt \rule{0pt}{12pt}\gen 12 \hskip .3em \id \hskip 5pt}\right)
\op 
\span\left(
\frac{\gen 12\hskip .2em \gen 12}{\gen 21\hskip .2em \twotwo}
\right)
\op
\span\left(\frac{\frac{\gen 12 \gen 12}{\rule{0pt}{9pt}\gen 21 \gen 21}}{
\hskip 5pt \rule{0pt}{12pt}\id \hskip .3em \gen 12\hskip 5pt}\right),
\]
where again $\raisebox{-2pt}{\twotwo} \in \S_0(2,2)$  is an arbitrary
element, modulo relations $R(2,3)$ that identify
the up-$(1,1;2)$-reducible generator of the second summand
with 
\[
\frac{\gen 12 \hskip .2em \gen 12}{\dvojiteypsilon \hskip .2em \gen 21} \in
\span\left(\frac{\gen 12 \hskip .2em \gen 12}{\twotwo \hskip .2em \gen 21}
\right) 
\]
and the up-$(1,1;2)$-reducible generator of the fourth summand with
\[
\frac{\gen 12 \hskip .2em \gen 12}{\gen 21 \hskip .2em \dvojiteypsilon} \in
\span\left(\frac{\gen 12 \hskip .2em \gen 12}{\gen 21  \hskip .2em \twotwo}
\right). 
\]
Finally, $\S_2(2,3)$ is the quotient of
\[
\span\left(\frac{\gen 12 \hskip 18pt \gen 12}{
\frac{\gen 12 \gen 12}{\rule{0pt}{9pt}\gen 21 \gen 21} 
\hskip 3pt \frac{\textstyle \id \id}{\rule{0pt}{9pt}\gen 21}
}\right)
\op
\span\left(
\frac{\hskip 3pt\raisebox{-2pt}%
{\onethree \hskip 3pt \onethree}}{\gen 21 \gen 21 \gen 21}
\right)
\op
\span\left(\frac{\gen 12 \hskip 18pt\gen 12}{
\frac{\textstyle\id \id}{\rule{0pt}{9pt}\gen 21}
\hskip 3pt \frac{\gen 12 \gen 12}{\rule{0pt}{9pt}\gen 21 \gen 21}
}\right)
\]
modulo $R(2,3)$ identifying the up-$(2,1;1,1)$-reducible element
\[
\frac{\ZbbZb \hskip 3pt \ZbbZb}{\gen 21 \gen 21 \gen 21}
\in
\span\left(
\frac{\hskip 3pt\raisebox{-2pt}%
{\onethree \hskip 3pt \onethree}}{\gen 21 \gen 21 \gen 21}
\right)
\]
with the generator of the first
summand, and the up-$(1,2;1,1)$-reducible element
\[
\frac{\bZbbZ \hskip 3pt \bZbbZ}{\gen 21 \gen 21 \gen 21}
\in
\span\left(
\frac{\hskip 3pt\raisebox{-2pt}%
{\onethree \hskip 3pt \onethree}}{\gen 21 \gen 21 \gen 21}
\right)
\]
with the generator of the last summand.
}
\end{example}

Example~\ref{kourim_vodni_dymku} shows that
presentation~(\ref{uz_je_hodne_pozde}) is not economical. Moreover, we
do not need to delve into the structure of $\S_0(m,n)$ because we
already know that this piece of $\S(m,n)$, isomorphic to the free
non-$\Sigma$ \hPROP\ $\uGamma_{\frac 12}(\uXi)$, is $\pa_0$-acyclic,
see Remarks~\ref{kasparek}, ~\ref{kasparek1} and
Theorem~\ref{jsem_prepracovany}. So we will work with the {\em
reduced form\/} of presentation~(\ref{uz_je_hodne_pozde}):
\begin{equation}
\label{uz_je_hodne_pozde-reduced}
\sfS(m,n) 
=  \S_0(m,n) \oplus  
\bigoplus_N \sfS(a_1,u) \ot \cdots \ot \sfS(a_v,u)
\ot
\sfS(v,b_1) \ot \cdots \ot \sfS(v,b_u)
/Q(m,n)
\end{equation}
where
\begin{equation}
\label{jsem_v_Bukuresti!}
N := M \cap \{v \geq 2,\ u\geq 2\} =\{
2 \leq v \leq m,\ 2 \leq u  \leq n, \mbox { and }
\textstyle\sum_1^v a_i =m,\ \sum_1^u b_j =n
\}
\end{equation}
and $Q(m,n) \subset R(m,n)$ is the span of $(\Rada c1s;\Rada
d1t)$-relations with $s,t \geq 2$. 

\begin{example}
{\rm\
We have the following reduced representations:
\begin{eqnarray}
\label{J}
\S(2,2) &=& \S_0(2,2) \op \span\left(
    \frac{\jednadva \hskip .2em \jednadva}{\dvajedna \hskip .2em\dvajedna}
    \right) \ \mbox { and}
\\
\nonumber 
\S_1(2,3) &=& \span\left(
\frac{\gen 12 \hskip .2em \gen 12}{\twotwo \hskip .2em \gen 21}
\right) 
\op 
\span\left(
\frac{\gen 12\hskip .2em \gen 12}{\gen 21\hskip .2em \twotwo}
\right), \ \mbox { where }\twotwo \in \sfS_0(2,2).
\end{eqnarray}
The reduced presentation of $\S_2(2,3)$ is the same as the unreduced
one given in Example~\ref{kourim_vodni_dymku}. We conclude that
\[
\S(2,3) \cong \S_0(2,3) \op \S_0(2,2) \op \S_0(2,2) \op \S_0(1,3)^{\ot
2}. 
\]
This, by the way, immediately implies the $\pa_0$-acyclicity of $\S(2,3)$.
}
\end{example}

\section{Acyclicity of the space of special elements}
\label{ac}

The proof of the $\pa_0$-acyclicity, in positive dimensions, of
$\S(m,n)$ is given by induction on $K := m \cdot n$. The acyclicity is
trivial for $K  \leq 2$. Indeed, there are only three spaces to
consider, namely $\sfS(1,1) = \span(\id)$, 
$\sfS(1,2) = \span(\gen12)$ and $\sfS(2,1) =\span(\gen21)$. All these
spaces are concentrated in degree zero and have trivial differential. 

The acyclicity is, in fact, obvious also for $K=3$ because $\S(1,3) =
\S_0(1,3)$ and $\S(3,1) = \S_0(3,1)$ coincide with their tree
parts. For $K=4$ we have two `easy' cases, $\S(4,1) = \S_0(4,1)$ and
$\S(1,4) = \S_0(1,4)$, while the acyclicity of $\S(2,2)$ follows from
presentation~(\ref{J}). 

Suppose we have proved the $\pa_0$-acyclicity of all $\sfS(k,l)$ with
$k \cdot l < K$. Let us express the reduced
presentation~(\ref{uz_je_hodne_pozde-reduced}) as the short exact
sequence
\begin{equation}
\label{zavody_probihaji}
0 \lra Q(m,n) \lra L(m,n) \oplus \sfS_0(m,n) \lra \sfS(m,n) \lra 0,
\end{equation}
where we denoted
\[
L(m,n) := \bigoplus_N \sfS(a_1,u) \ot \cdots \ot \sfS(a_v,u)
\ot \sfS(v,b_1) \ot \cdots \ot \sfS(v,b_u)
\]
with $N$ defined in~(\ref{jsem_v_Bukuresti!}).  It follows from the
K\"unneth formula and induction that $L(m,n)$ is $\pa_0$-acyclic while
the acyclicity of $\sfS_0(m,n) \cong \uGamma_{\frac 12}(\uXi)$ was
established in Remark~\ref{kasparek}.  Short exact
sequence~(\ref{zavody_probihaji}) then implies that it is in fact
enough to prove that the space of relations $Q(m,n)$ is
$\pa_0$-acyclic for any $m,n \geq 1$. This would clearly follow from
the following claim.

\begin{claim}
\label{25}
For any $w \in L(m,n)$ such that $\pa_0 (w) \in Q(m,n)$, 
there exists $z \in Q(m,n)$ such that $\pa_0(z) = \pa_0(w)$. 
\end{claim}

\begin{proof}
It follows from the nature of relations in the reduced
presentation~(\ref{uz_je_hodne_pozde-reduced}) that 
\begin{equation}
\label{H}
\pa_0(w) = \sum_{\cd} u_\uparrow^{\cd} - u_\downarrow^{\cd},
\end{equation}
where the summation runs over all $\cd = (\Rada c1s;\Rada d1t)$ with
$s,t \geq 2$, and $u_\uparrow^{\cd}$ (resp. $u_\downarrow^{\cd}$) is
an \cd-up (resp.~down) reducible element such that $u_\uparrow^{\cd} -
u_\downarrow^{\cd} \in Q(m,n)$. The idea of the proof is to show that
there exists, for each \cd, some \cd-up-reducible $z_\uparrow^{\cd}$
and some \cd-down-reducible $z_\downarrow^{\cd}$ such that $z^{\cd} :
= z_\uparrow^{\cd} - z_\downarrow^{\cd}$ belongs to $Q(m,n)$ and
\begin{equation}
\label{P}
u_\uparrow^{\cd} - u_\downarrow^{\cd} = \pa_0 (z^{\cd}).
\end{equation}
Then $z := \sum_\cd z^{\cd}$ will certainly fulfill $\pa_0(z) =
\pa_0(w)$. We will distinguish five types of $\cd$'s. 
The first four types are easy to handle; the last type is more
intricate.

{\it Type 1: All $\Rada d1t$ are $\geq 2$ and all $\Rada c1s$ are
arbitrary.} In this
case $u_\uparrow^\cd$ is of the form
\begin{equation}
\label{F}
\frac{A_1 \cdots A_t}{B_1 \cdots B_u} 
\end{equation}  
with 
\[
A_i = \frac{A_{i1} \cdots A_{id_i}}{C_{i1} \cdots C_{is}}
\]
as in~(\ref{U1}).  It follows from the definition that {\em $\pa_0$
cannot create $(k,l)$-fractions with $k,l \geq 2$.} Therefore a
monomial as in~(\ref{F}) may occur among monomials forming $\pa_0(y)$
for some monomial $y$ if and only if $y$ itself is of the above
form. Let $z^{\cd}_\uparrow$ be the sum of all monomials in $w$ whose
$\pa_0$-boundary nontrivially contributes to $u_\uparrow^\cd$. Let
$z^{\cd}_\downarrow$ be the corresponding \cd-down-reducible
element. Then clearly $u_\downarrow^\cd = \pa_0(z_\downarrow^\cd)$
and~(\ref{P}) is satisfied with $z^\cd := z_\uparrow^\cd -
z^\cd_\downarrow$ constructed above. In this way, we may eliminate all
$\cd$'s of Type~1 from~(\ref{H}).

{\it Type 2: All $\Rada c1t  = 1$ and all $\Rada d1s$ are arbitrary.}
In this case $u_\uparrow^\cd$ is of the form
\begin{equation}
\label{Z}
\frac{A_1 \cdots A_t}{B_1 \cdots B_u}
\end{equation}
with
\begin{equation}
\label{B}
A_i = \frac{A_{i1} \cdots A_{id_i}}{C_{i1} \cdots C_{is}},
\end{equation}
where $C_{ij} \in \S(d_i,1)$, for $1 \leq i \leq t$, $1 \leq j \leq
s$. When $d_i = 1$,
\[
A_i = A_{i1} \circ (C_{i1} \ot \cdots \ot C_{is}),
\]
where $C_{ij} \in \S(1,1)$ must be a scalar multiple of $\id$. We
observe that element~(\ref{Z}) is $\cd$-up reducible if and only if
$A_i$ is as in~(\ref{B}) where $d_i \geq 2$; if $d_i = 1$ then $A_i$ may
be arbitrary. We conclude, as in the previous case, that a monomial of
the above form may occur in $\pa_0(y)$ if and only if $y$ itself is of
the above form. Therefore, by the same argument, we may eliminate
$\cd$'s of Type~2 from~(\ref{H}).

{\it Type 3: All $\Rada c1t \geq 2$ and $\Rada d1s$ are arbitrary.} This
case is dual to Type~1.

{\it Type 4: All $\Rada d1t  = 1$ and $\Rada c1s$ are arbitrary.} This case is
dual to Type~2.

{\it Type 5: The remaining case.} This means that $1 \in \{\Rada
c1s\}$ but there exist some $c_j \geq 2$, and $1 \in \{\Rada d1t\}$ but
there exists some $d_i \geq 2$. This is the most intricate case, because
it may happen that, for some monomial $y$, $\pa_0(y)$ contains a \cd- (up-
or down-) reducible piece although $y$ itself is not
$\cd$-reducible. For instance, let 
\[
y := \frac{
           \frac {\raisebox{2pt}{ \zeroatwo {} \hskip -.65em\zeroatwo {}} }
           {\twoatwo {} \hskip .25em\twoaone {}} 
{
\unitlength=.4pt
\begin{picture}(98.00,60.00)(0.00,25.00)
\put(40.00,40.00){\line(1,-1){40.00}}
\put(40.00,40.00){\line(-1,-1){40.00}}
\put(40.00,60.00){\line(0,-1){60.00}}
\end{picture}}
          }
          {\twoazero {} \hskip .85em \twoazero {} \hskip .85em\twoazero {}}.
\]
Then $\pa_0(y)$ contains a $(2,1;2,1)$-up-reducible piece
\[
\frac{
           \frac {\raisebox{2pt}{ \zeroatwo {} \hskip -.65em\zeroatwo {}} }
           {\twoatwo {} \hskip .25em\twoaone {}} 
{
\unitlength=.4pt
\begin{picture}(98.00,60.00)(0.00,25.00)
\bezier{50}(20.00,20.00)(30.00,10.00)(40.00,0.00)
\put(40.00,60.00){\line(0,-1){20.00}}
\put(40.00,40.00){\line(1,-1){40.00}}
\put(40.00,40.00){\line(-1,-1){40.00}}
\end{picture}}
}
          {\twoazero {} \hskip .85em \twoazero {} \hskip .85em\twoazero {}}
\]
though $y$ itself is not $(2,1;2,1)$-up-reducible.
Nevertheless, we see that $\pa_0(y)$ contains also
\[
\frac{
           \frac {\raisebox{2pt}{ \zeroatwo {} \hskip -.65em\zeroatwo {}} }
           {\twoatwo {} \hskip .25em\twoaone {}} 
{
\unitlength=.4pt
\begin{picture}(98.00,60.00)(0.00,25.00)
\bezier{50}(55.00,25.00)(40.00,10.00)(30.00,0.00)
\put(40.00,60.00){\line(0,-1){20.00}}
\put(40.00,40.00){\line(1,-1){40.00}}
\put(40.00,40.00){\line(-1,-1){40.00}}
\end{picture}}
}
          {\twoazero {} \hskip .85em \twoazero {} \hskip .85em\twoazero {}}
\]
which is \underline{not} reducible. This is in fact a general
phenomenon, that is, if $\pa_0(y)$ contains a \cd-up-reducible piece
and if $y$ is not \cd-up-reducible, then $\pa_0(y)$ contains also an
irreducible piece. Therefore such $y$ cannot occur among monomials
forming up $w$ in  Claim~\ref{25}. We conclude that $y$ must also be
\cd-up-reducible and eliminate it from~(\ref{H}) as in the previous
cases. Down-reducible pieces can be handled similarly.
This finishes our proof of Claim~\ref{25}.%
\end{proof}

\section{Some generalities on minimal models}
\label{po_zavodech_v_Hosine}

In this section we discuss properties of minimal models of \PROP{s}.
We will see that minimal models of \PROP{s} do not behave as nicely as
for example minimal models of simply connected commutative associative
algebras.  We will start with an example of a \PROP\ that does not
admit a minimal model.  Even when a minimal model of a given \PROP\
exists, we are not able to prove that it is unique up to isomorphism,
although we will show that it is still unique in a weaker sense.
These pathologies of minimal models for \PROP{s} are related to the
absence of a suitable filtration required by various inductive
procedures used in the ``standard'' theory of minimal models.

In this section we focus on minimal models of \PROP{s} that are
concentrated in (homological) degree $0$. This generality would be enough for
the purposes of this paper. Observe that even these very special
\PROP{s} need not have minimal models. 
An example is provided by the \PROP\
\[
\sfX := \Gamma(u,v,w)/(u\circ v = w, v\circ w = u, w\circ u = v),
\]
where $u$, $v$ and $w$ are degree $0$ generators of biarity $(2,2)$.
Before we show that $\sfX$ indeed does not admit a minimal model,
observe that a (non-negatively graded) minimal model of 
an arbitrary \PROP\ concentrated in degree $0$ is always of the form
\[
\sfM = (\Gamma(E),\pa),
\]
where $E = \bigoplus_{i \geq 0} E_i$ with $E_i := \{e \in E;\ \deg(e) =
i\}$, and the differential $\pa$ satisfying, for any $n \geq 0$,
\begin{equation}
\label{77}
\pa(E_n) \subset \Gamma(E_{< n}),\ E_{< n} := \textstyle\bigoplus_{i <n}E_i.
\end{equation}
This means that $\sfM$ is {\em special cofibrant\/} in the sense
of~\cite[Definition~17]{markl:ha}. 

Free \PROP{s} $\Gamma(E)$ are canonically augmented, with the
augmentation defined by $\epsilon(E) = 0$.  This augmentation induces
an augmentation of the homology of minimal dg-\PROP{s}, therefore all
\PROP{s} with trivial differential which admit a minimal model are
augmented. The contrary is not true, as shown by the example of the
\PROP\ $\sfX$ above with the augmentation given by $\epsilon(u) =
\epsilon(v) = \epsilon(w) := 0$.

Indeed, assume that $\sfX$ has a minimal model $\rho : (\Gamma(E),\pa) \to
(\sfX,0)$. The map $\rho$ induces the isomorphism
\[
H_0(\rho) : H_0(\Gamma(E),\pa) = \Gamma(E_0)/(\pa (E_1))
\stackrel{\cong}{\longrightarrow} \sfX.
\]
Since $\sfX = \bfk \oplus \sfX(2,2)$, $E_0 = E_0(2,2)$ and the above
map is obviously an isomorphism of augmented \PROP{s}.
Therefore $H_0(\rho)$ induces an isomorphism of the spaces of
indecomposables. While it follows from the minimality of $\pa$
that $Q(\Gamma(E_0)/(\pa (E_1))) \cong E_0$, clearly $Q(\sfX) = 0$,
from which we conclude that $E_0 = 0$, which is impossible.

Although we are not able to prove that minimal models are unique up to
isomorphism, the following theorem shows
that they are still well-defined objects of a certain derived
category. Namely, let $\catho-dgPROP$ be the localization of the
category $\catdgPROP$ of differential non-negatively graded \PROP{s}
by homology isomorphisms.

\begin{proposition}
\label{pozitri_Hosin}
Let $\sfA$ be a \PROP\ concentrated in degree $0$. Then its minimal
model (if exists), considered as an object of the localized category
$\catho-dgPROP$, is unique up to isomorphism.
\end{proposition}

\noindent
{\bf Proof.}  The proposition would clearly be implied by the  
following statement.  Let $\alpha : \sfM' \to
\sfA$ and $\beta : \sfM'' \to \sfA$ be two minimal models of
$\sfA$. Then there exists a homomorphism $h : \sfM' \to \sfM''$ such
that the diagram
\[
\Vtriangle{\sfM'}{\sfM''}{\sfA}h{\alpha}{\beta}
\]
commutes.  Such a map $h$ can be constructed by induction.  Assume that $\sfM'
= (\Gamma(E),\pa')$, $\sfM'' = (\Gamma(F),\pa'')$ and let $h_0 :
\Gamma(E_0) \to \sfM''$ be an arbitrary lift in the diagram
\[
\Vtriangle{\Gamma(E_0)\ }{\sfM''}%
{\sfA}{h_0}{\alpha|_{\Gamma(E_0)}}{\beta}
\]
Suppose we have already constructed, for some $n \geq 1$, a
homomorphism
\[
h_{n-1}: \Gamma(E_{< n}) \to \sfM'',
\]
such that $\beta \circ h_{n-1} = \alpha|_{\Gamma(E_{< n})}$, where
$E_{< n}$ is as in~(\ref{77}).  Let us show that $h_{n-1}$ can be
extended into $h_n: \Gamma(E_{< n +1}) \to \sfM''$ with the similar
property.  To this end, 
fix a $\bfk$-linear basis ${\mathcal B}_n$ of $E_n$ and
observe that for each $e \in {\mathcal B}_n$ there
exists a solution $\omega_e \in \sfM''$ of the equation
\begin{equation}
\label{za_tri_dny_Hosin}
\pa'' (\omega_e) = h_{n-1}(\pa' e).
\end{equation}
This follows from the fact that $\pa'' h_{n-1}(\pa' e) = h_{n-1}(\pa
'\pa' e) = 0$ which means that the right-hand side
of~(\ref{za_tri_dny_Hosin}) is a $\pa''$-cycle, therefore $\omega_e$
exists by the acyclicity of $\sfM''$ in degree $n-1$. 
Define a linear map ${\underline r}_n :E_n \to \sfM''$ by ${\underline
r}_n(e) := \omega_e$, for $e \in {\mathcal B}_n$. Finally, define a linear
{\em equivariant\/} map $r_n :E_n \to \sfM''$ by
\[
r_n(f) := \sum_{\tau,\sigma} 
\frac 1{k!  \hskip .2em  l!} 
\hskip .2em \sigma^{-1} {\underline r}_n(\sigma f
\tau) \tau^{-1},
\]
where $f\in E_n$ is of biarity $(k,l)$ and the summation runs over all
$\sigma \in \Sigma_k$ and $\tau \in \Sigma_l$.  It is easy to verify
that the homomorphism $h_n : \Gamma(E_n) \to \sfM''$ determined by 
$h_n(f) := r_n(f)$ for $f \in E_n$, extends $h_{n-1}$,
and the induction goes on.

By modifying the proof of~\cite[Lemma~20]{markl:ha}, one may
generalize Proposition~\ref{pozitri_Hosin} to an arbitrary
non-negatively graded \PROP\ with trivial differential.

\begin{remark}
\label{99}
{\rm One usually proves that two minimal models connected by a
(co)homology isomorphisms are actually isomorphic. This is for
instance true for minimal models of 1-connected commutative
associative dg-algebras~\cite[Theorem~11.6(iv)]{lehmann:Asterisque77},
minimal models of connected dg-Lie
algebras~\cite[Theorem~II.4(9)]{tanre83} as well as for minimal models
of augmented operads~\cite[Proposition~3.120]{markl-shnider-stasheff:book}.
We do not know whether this isomorphism theorem is true also for
minimal models of \PROP{s}.
}
\end{remark}

However, `classical' isomorphism theorems can still be proved if one
imposes some additional assumptions on the type of minimal models
involved, as illustrated by Theorem~\ref{zitra_do_Hosina} below.  Let
us call a minimal model $(\Gamma(\Xi),\pa)$ of the bialgebra \PROP\
$\sfB$ {\em special\/} if the differential $\pa$ preserves the
subspace of special elements and if it is of the form $\pa = \pa_0 +
\papert$, where $\pa_0$ is as in~(\ref{Paal_v_Praze}) and $\papert$
raises the genus.

\begin{theorem}
\label{zitra_do_Hosina}
Let $\sfM' = (\Gamma(\Xi),\pa')$ and $\sfM'' = (\Gamma(\Xi),\pa'')$ be
two special minimal models for the bialgebra \PROP\ $\sfB$. Then
there exists an isomorphism $\phi: \sfM' \to \sfM''$ preserving the
space of special elements.
\end{theorem}

\noindent 
{\bf Proof.}  Let us construct inductively a map $\phi :
(\Gamma(\Xi),\pa') \to (\Gamma(\Xi),\pa'')$ of augmented \PROP{s}
which preserves the space of special elements and which is the `identity
modulo elements of higher genus.' By this we mean that $\phi|_\Xi =
\id_\Xi + \eta$, where the image of the linear map $\eta: \Xi \to \Gamma(\Xi)$
consists of special elements of positive genera.

The first step of the inductive construction is easy: we 
define $\phi_0 : \Gamma(\Xi_0)  \to
\Gamma(\Xi)$ by $\phi_0|_{\Xi_0} := \id_{\Xi_0}$.  
Suppose we have already constructed $\phi_{n-1} : \Gamma(\Xi_{< n})
\to \Gamma(\Xi)$. As in the proof of Proposition~\ref{pozitri_Hosin},
one must solve, for each element $e$ of a basis of $\Xi_n$, the equation
\begin{equation}
\label{111}
\pa'' \omega_e = \phi_{n-1} (\pa ' e).
\end{equation} 
The right-hand side is a $\pa''$-cycle, therefore a solution
$\omega_e$ exist by the acyclicity of $\sfM''$ in degree $n-1$.  But not
every solution is good for our purposes. Observe that the right-hand
side of~(\ref{111}) is of the form
\[
\phi_{n-1}(\pa_0 e + \papert' e) = \pa_0(e) +  \vartheta_e,
\]
where $\vartheta_e$ is a sum of special elements of positive genera.
We leave as an exercise to prove that the $\pa_0$-acyclicity of the
space of special elements implies, similarly as in the `na\"{\i}ve'
proof of Theorem~\ref{main} given in Section~\ref{zase_starosti_mam},
that in fact  there exists
$\omega_e$ of the form $\omega_e = e +
\eta_e$, where $\eta_e$ is a sum of special elements of
genus $> 0$. Therefore
\[
\phi_n(e) := \omega_e = e + \eta_e, \ e \in \Xi_n,
\]
defines an extension of $\phi_{n-1}$ which preserves special elements
and which is the identity modulo elements of higher genera.

The proof is concluded by showing that every endomorphism $\phi :
\Gamma(\Xi) \to \Gamma(\Xi)$ whose linear part is the identity and which
preserves the space of special elements is invertible.  We leave  
this statement as another exercise to the reader.%
\qed

\section{Proof of the main theorem and final reflections}
\label{proof_problems}

\begin{proof}[Proof of Theorem~\ref{main}]  
We already know from
Theorem~\ref{jedeme_do_Jinolic} that the collection $\sfS$ of special
elements is friendly, therefore the inductive construction described
in Section~\ref{zase_starosti_mam} gives a perturbation $\papert =
\pa_1 + \pa_2 + \pa_3 + \cdots$ \hskip 1mm such that
\[
\pa_g(\xi^m_n) \in \sfS_g(m,n), \mbox { for } g \geq 0.
\]
Equation~(\ref{L}) then immediately follows from Lemma~\ref{8} while
the fact that $\papert$ preserves the path grading follows from
Lemma~\ref{strasna_vedra}. 

It remains to prove that $({\sf M},\pa) = (\Gamma(\Xi),\pa_0 +
\papert)$ really forms a minimal model of~$\sfB$, that is, to
construct a homology isomorphism from $(\sfM,\pa)$ to $(\sfB,\pa
=0)$. To this end, consider the homomorphism
\[
\rho : (\Gamma(\Xi),\pa_0 + \papert) \to (\sfB,\pa = 0)
\]
defined, in
presentation~(\ref{pekne_jsem_si_ve_Zbraslavicich_poletal}), by
\[
\rho(\xi^1_2) := \jednadva, \
\rho(\xi^2_1) := \dvajedna,
\]
while $\rho$ is trivial on all remaining generators. It is clear that
$\rho$ is a well-defined map of dg-\PROP{s}.  The fact that $\rho$ is
a homology isomorphism follows from rather deep Corollary~27
of~\cite{mv}. An important assumption of this Corollary is that
$\papert$ preserves the path grading. This assumption guarantees that
the first spectral sequence of~\cite[Theorem~24]{mv} converges because
of the inequalities given in~\cite[Exercise~21]{mv} and recalled here
in~(\ref{netekla_voda}). The proof of Theorem~\ref{main} is finished.%
\end{proof}

\vskip 3mm
\noindent
{\bf Final reflections and problems.}  We observed that it is
extremely difficult to work with free \PROP{s}.  Fortunately, it turns
out that most of classical structures are defined over simpler objects
-- operads, \hPROP{s} or dioperads. In Remark~\ref{pocasi_nevychazi}
we indicated a definition of {\em special \PROP{s}\/} for
which only compositions given by `fractions' are allowed.

Let us denote by $\sB$ the special \PROP\ for bialgebras. It
clearly fulfills  $\sB(m,n) = {\bf k}$
for all $m,n \geq 1$ which means that bialgebras are 
the easiest objects defined over special \PROP{s} in the same sense in
which associative algebras are the easiest objects defined over
non-$\Sigma$-operads (recall that the non-$\Sigma$-operad 
$\underline{\Ass}$ for
associative algebras fulfills 
$\underline{\Ass}(n) = {\bf k}$ for all $n \geq 1$) and associative
commutative algebras are the easiest objects defined over
($\Sigma$-)operads (operad $\Comm$ fulfills $\Comm(n) = 
{\bf k}$ for all $n \geq 1$).

Let us close this paper by summarizing some open problems.

(1) Does there exist a sequence of convex polyhedra $B^m_n$
with the properties stated in Conjecture~\ref{v_patek_prileti_Jim}? 

(2) What can be said about the minimal model for the \PROP\ for
``honest'' Hopf algebras with an antipode?

(3) Explain why the Saneblidze-Umble diagonal occurs in our formulas
for $\pa$.

(4) Describe the generating function
\[
f(s,t) := \sum_{m,n \geq 1} \dim\S(m,n) s^m t^n
\]
for the space of special elements.

(5) Give a closed formula for the differential $\pa$ of the minimal model.

(6) Develop a theory of homotopy invariant versions of algebraic
objects over \PROP{s}, parallel to that of~\cite{markl:ha} for
algebras over operads. We expect that all main results
of~\cite{markl:ha} remain true also for \PROP{s}, though there might
be surprises and unexpected difficulties related to the combinatorial
explosion of \PROP{s}.

(7) What can be said about the uniqueness of the minimal model? Is the
minimal model of an augmented \PROP\ concentrated in degree $0$
unique up to isomorphism? If not, is at least a suitable
completition of the minimal model unique?

There is a preprint~\cite{saneblidze-umble:bialgebras} which might
contain answers to Problems~(1) and~(5).

\def\cprime{$'$}

\end{document}